\documentclass[12pt, english]{amsart}
\usepackage{amsmath}
\usepackage[toc,page]{appendix}
\usepackage[english]{babel}
\usepackage{mathrsfs}
\usepackage{amsmath,amssymb,mathtools}
\usepackage[all]{xy}
\usepackage{url}
\usepackage{smfthm}
\usepackage{extpfeil}
\usepackage{mathtools}
\usepackage{yhmath}
\usepackage{bm}
\usepackage{mathptmx}
\usepackage{extarrows}

\usepackage{geometry}
\usepackage[colorlinks=true,hyperindex, linkcolor=magenta, pagebackref=false, citecolor=cyan,pdfpagelabels]{hyperref}
\usepackage[capitalize]{cleveref}
\usepackage{geometry}
\newcommand{\scrmod}{\mathrm{SCRMod}}

\newcommand{\sfp}{\mathrm{sfp}}

\newcommand{\Fr}{\mathrm{Fr}}
\newcommand{\Mod}{\mathrm{Mod}}
\newcommand{\N}{\mathbf N}

\newcommand{\cp}{\mathrm{cp}}

\newcommand{\Lcal}{\mathcal{L}}
\def\drCl{\mathfrak{C}\mathfrak{l}}

\newcommand{\R}{\mathbf{R}}

\newcommand{\Scal}{\mathcal{S}}

\newcommand{\Hom}{\mathrm{Hom}}

\newcommand{\Zbb}{\mathbb{Z}}

\newcommand{\Cart}{\mathrm{Cart}}

\newcommand{\dg}{\mathrm{dg}}
\newcommand{\drm}{\mathrm{d}}

\newcommand{\scralg}{\mathrm{SCRAlg}}

\newcommand{\ccart}{\mathrm{ccart}}

\newcommand{\Bcal}{\mathcal{B}}
\newcommand\pr{\mathrm{pr}}

\newcommand{\ev}{\mathrm{ev}}
\newcommand{\id}{\mathrm{id}}

\newcommand{\f}{\mathbf f}

\newcommand{\an}{\mathrm{an}}

\newcommand{\Sym}{\mathrm{Sym}}

\newcommand{\fin}{\mathrm{fin}}
\newcommand{\Ker}{\mathrm{Ker}}

\newcommand{\CAlg}{\mathrm{CAlg}}
\newcommand{\heart}{\heartsuit}

\newcommand{\Pbb}{\mathbb{P}}
\newcommand{\Ocal}{\mathcal{O}}

\newcommand{\Spec}{\mathrm{Spec}}
\newcommand{\Fcal}{\mathcal{F}}
\newcommand{\Gcal}{\mathcal{G}}
\newcommand{\Set}{\mathrm{Set}}
\newcommand{\poly}{\mathcal{P}\mathrm{oly}}
\newcommand{\Acal}{\mathcal{A}}
\newcommand{\Alg}{\mathrm{Alg}}
\newcommand{\Cl}{\mathrm{Cl}}
\newcommand{\T}{\mathrm{T}}
\newcommand{\Pcal}{\mathcal{P}}

\newcommand{\Ffrak}{\mathfrak{F}}

\newcommand{\Hrm}{\mathrm{H}}
\newcommand{\bi}{\mathrm{bi}}

\newcommand{\Ccal}{\mathcal{C}}
\newcommand{\QCoh}{\mathrm{QCoh}}
\newcommand{\Dcal}{\mathcal{D}}
\newcommand{\op}{\mathrm{op}}
\newcommand{\D}{\mathrm{D}}
\newcommand{\Fun}{\mathrm{Fun}}
\newcommand{\h}{\mathrm{h}}

\newcommand{\Lbf}{\mathbf{L}}

\newcommand{\Proj}{\mathrm{Proj}}
\newcommand{\Hcal}{\mathcal{H}}
\newcommand{\Ani}{\mathrm{Ani}}

\newcommand{\Cat}{\mathrm{Cat}}

\newcommand{\Yfrak}{\mathfrak{Y}}
\newcommand{\Xfrak}{\mathfrak{X}}
\newcommand{\Fbb}{\mathbb{F}}
\newcommand{\Sp}{\mathrm{Sp}}
\newcommand{\lrm}{\mathrm{l}}
\newcommand{\Ebb}{\mathbb{E}}

\newcommand{\Lrm}{\mathrm{L}}
\newcommand{\Ecal}{\mathcal{E}}

\usepackage{amsmath,amssymb,tikz-cd}
\newcommand{\drpoly}{\mathfrak{P}\mathrm{oly}}
\newcommand{\Prm}{\mathrm{P}}

\newcommand{\For}{\mathrm{For}}
\newcommand{\dst}{\mathrm{dSt}}
\newcommand{\alg}{\mathrm{alg}}

\newtheorem{theorem}[subsubsection]{Theorem}
\newtheorem{proposition}[subsubsection]{Proposition}

\newtheorem{lemma}[subsubsection]{Lemma}

\newtheorem{corollary}[subsubsection]{Corollary}

\theoremstyle{definition}
\numberwithin{equation}{subsubsection}
\newtheorem{definition}[subsubsection]{Definition}

\newtheorem{remark}[subsubsection]{Remark}
\newtheorem{example}[subsubsection]{Example}
\newtheorem{construction}[subsubsection]{Construction}

\makeindex             

\begin{document}
\title{Generalized Clifford algebras, weighted polynomial laws, and adjunctions}
\author{Nguyen Xuan Bach}
\address{Education Zone, Hoa Lac Hi-tech Park, Km29, Thang Long Boulevard, Thach Hoa, Thach That, Ha Noi, Vietnam} 
\email{bachnx9@fe.edu.vn}
\begin{abstract}
In this paper, we study Clifford algebra construction from the perspective of adjunctions motivated by the general framework of Krashen and Lieblich. We introduce a category of weighted polynomial laws whose associated Clifford algebra functor is a left adjoint on this category. We introduce a notion of lifting weighted polynomial laws generalizing explicit constructions of Krashen and Lieblich. As another application, we construct derived Clifford algebras of derived homogeneous polynomial forms of arbitrary degree with noncommutative coefficients.
\end{abstract}
%
%
\maketitle
\tableofcontents 

\setcounter{section}{-1}

\section{Introduction}
To a quadratic form $Q:M\to R$ over an $R$-module $M$, one can associate its Clifford algebra $\Cl(M,R)$. This algebra is a classical mathematical object and is a fundamental tool in the theory of quadratic forms (see \cite{Knus_1991} for example). In the literature, the Clifford algebra is usually defined as the quotient of the tensor algebra $\T(M)$ by relations $$m\otimes m=Q(m).$$ Several extensions of this definition have been studied, in which one replaces quadratic forms by forms of higher degrees or families of forms (for example \cite{Roby_1968}, \cite{Haile_1984}, \cite{Pappa_2000}, \cite{Chapman_2015},  \cite{Kuo_2011}, \cite{Kul_2003}, \cite{HH_2007}). As is customary in algebra, their representations have also been considered. In some of these works (for example \cite{Chapman_2015}), a projective variety is introduced so that one can move from studying representations of a Clifford algebra to studying vector bundles satisfying certain conditions on this variety. This method can be traced back to \cite{BHU_1991} and \cite{Bergh_1987}; later, these vector bundles became known as Ulrich bundles. 

D. Krashen and M. Lieblich, in their preprint \cite{KL_2015}, are the first ones who take a very general point of view to unify these approaches. More precisely, for a 
morphism $\phi:X\to Y$ of schemes over a base scheme $S$, they consider the \textit{Clifford functor}
$$\Ccal\Fcal_\phi:\mathrm{LFAlg}/\Ocal_S \to \Set$$
defined by the formula $\Ccal\Fcal(\Bcal)=\Hom_{\Alg/\Ocal_Y}(\phi_*\Ocal_X,\Bcal|_Y),$ where $\mathrm{LFAlg}/\Ocal_S$ is the category of locally free sheaves of algebras over $S$. The algebra that corepresents $\Ccal\Fcal_\phi$ is called the Clifford algebra of $\phi$. For example, the classical Clifford algebra of a quadratic form is the Clifford algebra of a degree $2$ covering $\phi:X\to \Pbb^{n}$ associated to this quadratic form. Moreover, their explicit constructions considered as instances of the above general framework cover all previous generalized Clifford algebras such as the ones considered in \cite{Chapman_2015}.

Through their definition of the Clifford functor, one can see the idea of interpreting the Clifford algebra as some kind of adjunction. In the context of affine morphisms, there is equivalence between the category
affine morphisms to $Y$ and the category of quasi-coherent sheaves over $Y$
under which $\phi:X\to Y$ corresponds to $\phi_*\Ocal_X$. 
We then want to interpret Krashen and Lieblich's Clifford algebra construction as a left adjoint to the restriction functor $(-)|_Y$ between the category of quasi-coherent algebras over $S$ and the category of quasi-coherent algebras over $Y$, completely discarding the morphism $\phi:X\to Y$ in their construction. 

The goal of  this paper is to make this idea work. Our intuition is that this left adjoint exists because of the existence of an underlying adjunction at the level 
of modules. This allows us to work in the setting of monoidal stable $\infty$-categories and derived $\infty$-categories. Statements about quasi-coherent sheaves are then delivered by considering the hearts of canonical t-structures.

The main new definitions in our paper are a category of weighted polynomial laws over any scheme $S$ and its associated Clifford algebra functor $\Cl: \poly^d(S)\to \Alg_S$ the category of quasi-coherent sheaves of $\Ocal_S$-algebras. Roughly speaking, these categories consist of a polynomial law with coefficients in a module and finite number of polynomial laws with coefficients in associative algebras. A significant feature is that our Clifford functor is a left adjoint and our construction is in terms of colimits, rather than killing ideals. To our knowledge, this emphasis is new in the world of Clifford algebras, and it is quite convenient when dealing with functors like left adjoints, that commute with colimits. 

For any morphism of schemes $f:Y\to S$, there exists a pullback functor $f^*:\poly^d(S)\to \poly^d(Y).$
It is not known whether this pullback functor of weighted polynomial laws admits a left adjoint. The problem on the existence of such left adjoints gives rise to a notion of lifting weighted polynomial laws from $S$ to $Y$ along $f$ and we prove a theorem to provide such lifts. The main property of these lifts is that the left adjoint of $f^*$ maps $\Cl(\Psi_Y)$ to $\Cl(\Psi_S)$ where $\Psi_Y$ is a lift of $\Psi_S$. From this point of view, the connection between generalized Clifford algebras and algebraic geometry occurs when $\Cl(\Psi_Y)$ is commutative as one can instead work with its relative spectrum. 

In an other direction of application of the fact that our Clifford algebra functor is a left adjoint, we give the construction of a derived Clifford algebra of a derived homogeneous polynomial forms of any degree with noncommutative coefficients. Note that a derived version of the Clifford algebra for derived quadratic forms has already appeared in \cite{Vezz_2016}. 

Let us give a summary of the paper. In Section 1, we set up an adjunction between
the categories of algebras. More precisely, we provide some conditions under which
$\Alg(F) : \Alg(\Ccal) \to \Alg(\Dcal)$ admits a left adjoint for a monoidal functor $F : \Ccal\to \Dcal$ between
monoidal $\infty$-categories (Corollary \ref{cor_adjoint}). Then, we apply this adjunction to the (derived)
pullback functor between quasi-coherent derived $\infty$-categories to prove the existence of
the left adjoint of the pullback functor in Theorem \ref{ad_module}.

In Section 2, we introduce the category of weighted polynomial laws (Definition \ref{def:w_poly}) and the Clifford algebra functor (Theorem \ref{thm:Cliff}). We show a method to generate lifts (Theorem \ref{thm:main_compute}), and for the case of projective bundle $p:Y=\Pbb(\Fcal)\to S$ of a vector bundle $\Fcal$ over $S$, we obtain a lift of a specific weighted polynomial $\Psi$ over the base $S$ that we denote by $\Psi|_{\Ocal(-1)}$ (Construction \ref{con:restrict}). Then we want to be more explicit and make a precise connection with Krashen and Lieblich's explicit construction by giving the relation defining the relative spectrum $\phi:X=\Spec(\Psi|_{\Ocal(-1)})\to Y$. When $S$ is affine and $\Fcal$ is free, it coincides with their explicit construction given in \cite[Theorem A.1.4]{KL_2015}.

In Section 3, we give the construction of derived Clifford algebra associated to derived homogeneous polynomial laws and prove some basic properties. \\
\noindent
\textbf{Acknowledgements.} As the paper is part of my Phd thesis, I would like to express my gratitude to my thesis directors Baptiste Calm\`es and Anne Qu\'eguiner-Mathieu for their kindness and big support to me during the three years of my PhD study. Baptiste’s program of formulating the Clifford algebra construction at this level of generality, which they introduced to me, has been invaluable and I have learnt a lot of lessons from my supervisors during these years. At LAGA and Sorbonne Paris North University, I was provided with a comfortable environment to carry out the project, so I owe them a debt of gratitude. Thanks to Yonatan Harpaz for his suggestions of some important points in Section 1. I sincerely thank Daniel Krashen and Ben Williams who gave me their valuable opinions. I warmly thank Fr\'ed\'eric D\'eglise, Bruno Kahn, and Bruno Vallette for their kind concern regarding this project. I thank Trieu Thu Ha, Pham Ngo Thanh Dat, and Dao Quang Duc for mathematical discussions at the beginning of this project, and Hoang Manh Truong for teaching me some basics of higher category theory.
\section{Constructions of left adjoints}
\subsection{Adjunctions between symmetric monoidal $\infty$-categories}
We start by setting up an adjunction between monoidal $\infty$-categories. Recall that a monoidal $\infty$-category $\Ccal$ is \textit{closed} if for all $x\in \Ccal$, the functor $x\otimes_\Ccal -$ has a right adjoint (see \cite[Definition 4.1.1.15]{HA}). These right adjoints can be upgraded to a bifunctor $[-,-]:\Ccal^{\op}\times \Ccal\to \Ccal$. For convenience, denote by $\D_\Ccal$ the functor $[-,\id_\Ccal]:\Ccal \to \Ccal^{\op}.$ Let $F:\Ccal\to\Dcal$ be a monoidal functor between closed monoidal $\infty$-categories. We have the folllowing natural transformation 
$$(F\circ [-,-]_\Ccal\Rightarrow [-,-]_\Dcal\circ (F \times F)):\Ccal\times \Ccal\to \Dcal.$$ In particular, we have a natural transformation $$(F^{\op} \circ [-,\id_\Ccal]\Rightarrow [-,\id_\Dcal]\circ F):\Ccal\to\Dcal^{\op}.$$ We denote this transformaton by $\theta_F$. We say that $F$ 
\textit{commutes with duals} if the canonical transformation $$\theta_F:F^{\op} \circ \D_\Ccal\Rightarrow \D_\Dcal\circ F$$ is an isomorphism.
\begin{lemma}\label{module_adjoint}
Let $F:\Ccal \to \Dcal$ be a monoidal functor between closed monoidal $\infty$-categories with a right adjoint $G$. Assume that
\begin{enumerate}
\item The categories $\Ccal,\Dcal$ are compactly generated;
    
\item $\Ccal^{\omega}, \Dcal^{\omega}$ are closed monoidal $\infty$-subcategories of $\Ccal,\Dcal$ and the transformations $\id_{\Ccal^{\omega}}\Rightarrow\D_{\Ccal^{\omega}}^{\op}\D_{\Ccal^{\omega}}$ and $\id_{\Dcal^{\omega}}\Rightarrow \D_{\Dcal^{\omega}}^{\op}\D_{\Dcal^{\omega}}$ are isomorphisms;

\item $F$ maps $\Ccal^{\omega}$ into $\Dcal^{\omega}$ (or equivalently $G$ preserves filtered colimits) and $G$ maps $\Dcal^{\omega}$ into $\Ccal^\omega$;
\item And $F|_{\Ccal^{\omega}}:\Ccal^{\omega}\to \Dcal^{\omega}$ commutes with duals.
\end{enumerate}
Then $F$ admits a left adjoint and this left adjoint is given on $\Ccal^{\omega}$ by $\D_\Ccal^{\op}G^{\op}\D_{\Dcal^\omega}$.     
\end{lemma}
\begin{proof}
We have an adjunction $(F|_{\Ccal^{\omega}}\dashv G|_{\Dcal^{\omega}}): \Ccal^{\omega} \leftrightarrows\Dcal^{\omega}$. Then $F|_{\Ccal^{\omega}}$ admits a left adjoint $H:\Dcal^{\omega}\to \Ccal^{\omega}$. 
Indeed, we have the following diagram of functors
$$\xymatrix@R+1pc@C+1pc{
\Ccal^\omega \ar[rr]<2pt>^*{\D_{\Ccal^\omega}} \ar[d]<2pt>^*{F} && (\Ccal^\omega)^{\op} \ar[d]<2pt>^*{F^{\op}} \ar[ll]<2pt>^*{\D_{\Ccal^\omega}^{\op}}\\
\Dcal^\omega \ar[u]<2pt>^*{G}  \ar[rr]<2pt>^*{\D_{\Dcal^\omega}} && (\Dcal^\omega)^{\op} \ar[ll]<2pt>^*{\D_{\Dcal^\omega}^{\op}} \ar[u]<2pt>^*{G^{\op}}.
}$$
Noting that by assumptions $G^{\op}$ is a left adjoint of $F^{\op}$, one has that ${\D_{\Ccal^\omega}}^{\op}\circ G^{\op} \circ \D_{\Dcal^\omega}$ is a left adjoint of $\D_{\Dcal^\omega}^{\op}\circ F^{\op}\circ \D_{\Ccal^\omega}$ since $\D_{\Ccal^\omega},\D_{\Dcal^\omega}$ are equivalences. Moreover,
$$\D^{\op}_{\Dcal^\omega}\circ F^{\op}\circ \D_{\Ccal^\omega} \simeq \D_{\Dcal^\omega}^{\op}\circ \D_{\Dcal^\omega}\circ F\simeq F.$$
So we can put $H=\D_{\Ccal^\omega}^{\op}\circ G^{\op}\circ \D_\Dcal.$
The conclusion of the lemma is equivalent to showing that $\h^x\circ F$ is corepresentable for any $x\in \Dcal$. As $\Dcal$ is generated by $\Dcal^{\omega}$ under colimts, it suffices to show that $\h^x \circ F$ is corepresentable for any compact object $x\in \Dcal$. Pick any $y\in \Ccal$. There exists a filtered diagram $p:K\to \Ccal^{\omega}$ such that $\mathrm{colim} p\cong y.$ Then 
$$
\begin{matrix*}[l]
\h^{Hx}(y) &\simeq &\h^{Hx}(\mathrm{colim} p)\\
&\simeq &\mathrm{colim}(\h^{Hx}\circ p)&(Hx \text{ is compact})\\
&\simeq& \mathrm{colim}(h^x \circ F \circ p)& (F|_{\Ccal^{\omega}}\text{ is a right adjoint of }H)\\  
&\simeq &\h^x(\mathrm{colim}(F\circ p))&(x\text{ is compact})\\
&\simeq& \h^x(F(\mathrm{colim}(p)))&(F \text{ is a left adjoint})\\
&\simeq& \h^x(Fy),
\end{matrix*}
$$
which says that $\h^x \circ F$ is corepresentable. This completes the proof.  
\end{proof}
For a monoidal $\infty$-category $\Ccal$, we denote by $\Alg(\Ccal)$ the $\infty$-category of algebra objects in $\Ccal$ (see \cite[Notation 4.1.1.8]{HA} for concrete definitions), and if $F:\Ccal\to \Dcal$ is a lax monoidal functor between monoidal $\infty$-categories, we denote by $$\Alg(F):\Alg(\Ccal)\to \Alg(\Dcal)$$ the induced functor on $\infty$-categories of algebras. If $\Ccal$ is a monoidal ordinary category, then $\N(\Ccal)$ is a monoidal $\infty$-category and $\Alg(\N(\Ccal))\simeq \N(\Alg(\Ccal))$ where $\Alg(\Ccal)$ is the ordinary category of algebra objects $\Ccal$ (see \cite[Remark 1.1.16]{DAGII}).

\subsection{Adjunctions between $\infty$-categories of algebras}
Recall that by {\cite[Proposition 1.5.12]{DAGII}}, for any $\infty$-category $\Ccal$ which admits all small colimits equipped with a monoidal structure compatible with all small colimits, the forgetful functor $\theta:\Alg(\Ccal)\to \Ccal$ admits a left adjoint $\Fr$ which we consider as the free algebra functor. More explicitly, for every $A\in \Alg(\Ccal)$, there exists a simplicial object $A_\bullet$ of $\Alg(\Ccal)$ such that each $A_n$ is in the essential image of $\Fr$ and $|A_\bullet|\simeq A$. Moreover, $A_\bullet$ can be chosen so that $A_n\simeq T^{n+1}A$ where $T:=\Fr_\Ccal\circ \theta_\Ccal$. Also, $\Alg(\Ccal)$ admits all small colimits.
\begin{lemma}\label{main_adjoint}
Let $F:\Ccal \to \Dcal$ be a lax monoidal functor between monoidal $\infty$-categories which admit all small colimits. Assume that one has adjunction pairs $(H\dashv F)$ and and the monoidal structure of $\Dcal$ is compatible with all small colimits. Then $\Alg(F)$ admits a left adjoint.
\begin{proof}
Let $A$ be an arbitrary algebra in $\Alg(\Dcal)$. If $A$ is a free algebra, i.e., $\h^M\circ \theta$ is corepresented by $A$ for some object $M\in \Dcal$, then 
$$
\begin{matrix*}[l]
    \h^A\circ \Alg(F)& \simeq & \h^{M} \circ \theta_\Dcal \circ \Alg(F)\\
    & \simeq & \h^M \circ F\circ \theta_\Ccal\\
    &\simeq & \h^{H M} \circ \theta_\Ccal&(H \text{ is a left adjoint of F})\\
    & \simeq & \h^{\Fr H M}.
\end{matrix*}
$$

In general, as the monoidal structure of $\Dcal$ is compatible with all small colimits, any algebra $A\in \Alg(\Dcal)$ can be written as colimit of some diagram in the essential image of $\Fr_{\Dcal}$. Hence $\h^A\circ \Alg(F)$ is corepresentable and so $\Alg(F)$ admits a left adjoint. 
\end{proof}
\end{lemma}
\begin{corollary}\label{cor_adjoint}
Let $F:\Ccal \to \Dcal$ be a monoidal functor between closed monoidal $\infty$-categories which admits a right adjoint $G$. \begin{enumerate}
\item The categories $\Ccal,\Dcal$ are compactly generated;
    
\item $\Ccal^{\omega}, \Dcal^{\omega}$ are closed monoidal $\infty$-subcategories of $\Ccal,\Dcal$, and the transformations $\id_{\Ccal^{\omega}}\Rightarrow\D_{\Ccal^{\omega}}^{\op}\D_{\Ccal^{\omega}}$ and $\id_{\Dcal^{\omega}}\Rightarrow \D_{\Dcal^{\omega}}^{\op}\D_{\Dcal^{\omega}}$ are isomorphisms;

\item $F$ maps $\Ccal^{\omega}$ into $\Dcal^{\omega}$ (or equivalently $G$ preserves filtered colimits) and $G$ maps $\Dcal^{\omega}$ into $\Ccal^\omega$;

\item And $F|_{\Ccal^{\omega}}:\Ccal^{\omega}\to \Dcal^{\omega}$ commutes with duals. 
\end{enumerate}    
Then $\Alg(F)$ admits a left adjoint. 
\begin{proof}
The adjunction pair $(F\dashv G)$ together with the assumptions $(1)-(4)$ satisfy the assumptions of Lemma \ref{module_adjoint} so $F$ admits a left adjoint $H$. Then the adjunctions pairs $(H\dashv F\dashv G)$ imply that $\Alg(F)$ admits a left adjoint by Lemma \ref{main_adjoint}.        
\end{proof}
\end{corollary}
Recall that for an exact functor $F:\Ccal\to \Dcal$ between stable $\infty$-categories equipped with t-structures, we say that $F$ is \textit{right} (respectively \textit{left}) \textit{t-exact} if it maps $\Ccal_{\geq 0}$  (respectively $\Ccal_{\leq 0})$ to $\Dcal_{\geq 0}$ (respectively $\Dcal_{\leq 0}).$
\begin{remark}\label{monoidal_t}
\begin{enumerate}
    \item If $\Ccal$ is a closed monoidal and stable $\infty$-category then it is a stable monoidal $\infty$-category since for all $M\in \Ccal$, the functors $M\otimes -$ has a right adjoint so they are exact.
   \item For a stable monoidal $\infty$-category $\Ccal$ we have $\Alg(\Ccal_{\geq 0})\subseteq \Alg(\Ccal)$ is stable under colimits. Indeed, let $p:K^{\triangleright} \to \Alg(\Ccal)$ be a colimit diagram such that $p|_{K}$ factors through $\Alg(\Ccal_{\geq 0}).$ Then $\theta_{\Ccal}\circ p:K^{\triangleright}\to \Ccal$ is also a colimit diagram with $\theta_{\Ccal}\circ p|_{K}$ factors through $\Ccal_{\geq 0}$. As $\Ccal_{\geq 0}\subseteq \Ccal$ is stable under colimit, $\theta_{\Ccal}\circ p$ factors through $\Ccal_{\geq 0}$. It follows that $p$ factors through $\Alg(\Ccal_{\geq 0}).$
\end{enumerate}
\end{remark}
\begin{lemma}\label{exact}
Let assume the assumptions of Lemma \ref{main_adjoint}. Denote by $\tilde{H}$ the left adjoint of $\Alg(F).$ Assume further that $\Ccal$ and $\Dcal$ are stable monoidal $\infty$-categories equipped with compatible t-structures. If $F$ is left t-exact then $$\tilde{H}(\Alg(\Dcal_{\geq 0}))\subseteq \Alg(\Ccal_{\geq 0}).$$
\begin{proof}
Pick $A\in \Alg(\Dcal_{\geq 0})$. By \cite[Remark 3.5.6 ]{DAGII}, we can write $A$ as colimit of a simplicial object $A_\bullet$ where $A_n\simeq T^{n+1}A$ and $T=\Fr_\Dcal\theta_\Dcal$. Then $\tilde{H}A\simeq |\tilde{H}A_\bullet|$ and $\tilde{H}A_n\simeq \tilde{H}T^{n+1}A$. By Remark \ref{monoidal_t}, it suffices to show that $\tilde{H}T^{n+1}A\in \Alg(\Ccal_{\geq 0})$. As $TA\simeq\Fr_{\Dcal}\theta_{\Dcal}A$, we have $TA\in \Alg(\Dcal_{\geq 0})$ and so $T^n A\in \Alg(\Dcal_{\geq 0})$ for all $n$. We write
$$\tilde{H}T^{n+1}A\simeq\tilde{H}\Fr_{\Dcal}\theta_{\Dcal}T^{n}A\simeq \Fr_{\Ccal}H\theta_\Dcal T^n A.$$
Since $\theta_\Dcal T^n A \in \Dcal_{\geq 0}$, we have $H\theta_{\Dcal}T^n A\in \Ccal_{\geq 0}$ ($H$ is right t-exact) and hence $\Fr_\Dcal H\theta_\Dcal T^nA \in \Alg(\Ccal_{\geq 0}).$ 
\end{proof}
\end{lemma}
Let $\Ccal$ be a stable monoidal $\infty$-category equipped with a compatible t-structure. Then one can equip $\Ccal^{\heartsuit}$ a monoidal structure such that $\pi_0^{\Ccal}:=\tau_{\leq 0}|_{\Ccal_{\geq 0}}:\Ccal_{\geq 0}\to \Ccal^{\heartsuit}$ is a monoidal functor with a lax monoidal right adjoint $i_{\Ccal}:\Ccal^{\heartsuit}\to \Ccal_{\geq 0}$ (see \cite[Proposition 2.2.1.9]{HA}). This induces an adjunction $\Alg(\pi_0^{\Ccal})\dashv \Alg(i_\Ccal).$\begin{lemma}\label{lr_t-exact}
Let us be in the situation of Lemma \ref{exact}. If $F$ is both left t-exact and right t-exact then $\Alg(F|_{\Ccal^{\heartsuit}}):\Alg(\Ccal^{\heartsuit})\to \Alg(\Dcal^{\heartsuit})$ admits a left adjoint.
\begin{proof}
Noting that $F$ is left and right t-exact, we have the following commutative diagram    
$$
\xymatrix{
\Ccal^{\heartsuit} \ar[r]^{i_\Ccal} \ar[d]_{F|_{\Ccal^{\heartsuit}}} & \Ccal_{\geq 0} \ar[d]^{F|_{\Ccal_{\geq 0}}} \ar[r]& \Ccal \ar[d]^{F}\\
\Dcal^{\heartsuit} \ar[r]_{i_\Dcal} & \Dcal_{\geq 0} \ar[r]& \Dcal,
}
$$
which induces the commutative diagram:
$$
\xymatrix{
\Alg(\Ccal^{\heartsuit}) \ar[r]^{\Alg(i_\Ccal)} \ar[d]_{\Alg(F|_{\Ccal^{\heartsuit}})} & \Alg(\Ccal_{\geq 0}) \ar[d]^{\Alg(F|_{\Ccal_{\geq 0}})} \ar[r]& \Alg(\Ccal) \ar[d]^{\Alg(F)}\\
\Alg(\Dcal^{\heartsuit}) \ar[r]_{\Alg(i_\Dcal)} & \Alg(\Dcal_{\geq 0}) \ar[r]& \Alg(\Dcal).
}
$$
Denote the left adjoint of $\Alg(F)$ by $\tilde{H}$. By Lemma \ref{exact}, $\tilde{H}$ maps $\Alg(\Dcal_{\geq 0})$ to $\Alg(\Ccal_{\geq 0})$ so we have an adjunction pair $(\Alg(F|_{\Ccal_{\geq 0}})\dashv \tilde{H}|_{\Alg(\Dcal_{\geq 0})})$. Then $$\Alg(\pi_0^{\Ccal})\tilde{H}|_{\Alg(\Dcal_{\geq 0})}\Alg(i_\Dcal)$$ is a left adjoint of $\Alg(F|_{\Ccal^{\heartsuit}})$ as for $A \in \Alg(\Dcal^{\heartsuit})$ and $B\in \Alg(\Ccal^{\heartsuit})$, one has
$$
\begin{matrix*}
\Hom_{\Alg(\Ccal^{\heartsuit})}(\Alg(\pi_0^{\Ccal})\tilde{H}|_{\Alg(\Dcal_{\geq 0})}\Alg(i_\Dcal) A, B)
&\simeq & \Hom_{\Alg(\Ccal_{\geq 0})}(\tilde{H}|_{\Alg(\Dcal_{\geq 0})}\Alg(i_\Dcal)A,\Alg(i_\Ccal)B)\\
&\simeq  &\ \Hom_{\Alg(\Dcal_{\geq 0})}(\Alg(i_\Dcal)A,\Alg(F|_{\Ccal_{\geq 0}})\Alg(i_\Ccal)B)\\
&\simeq  &\ \Hom_{\Alg(\Dcal_{\geq 0})}(\Alg(i_\Dcal)A,\Alg(i_\Dcal)\Alg(F|_{\Ccal^{\heartsuit}})B)\\
&\simeq &\ \Hom_{\Alg(\Dcal^{\heartsuit})}(\Alg(\pi_0^{\Dcal})\Alg(i_\Dcal)A,\Alg(F|_{\Ccal^{\heartsuit}})B)\\
&\simeq 
&\ \Hom_{\Alg(\Dcal^{\heartsuit})}(A,\Alg(F|_{\Ccal^{\heartsuit}})B).
\end{matrix*}
$$
Here the first and the fourth equivalences are due to the adjunction pair $(\Alg(\pi_0)\dashv \Alg(i)),$ the third equivalence is due to our commutative diagrams, and the last equivalence is due to $\pi_0^{\Dcal}\circ i_\Dcal\simeq \id_{\Dcal^{\heartsuit}}.$
\end{proof}
\end{lemma} 
\subsection{Applications to quasi-coherent derived categories}
\begin{theorem}\label{ad_module}
Let $f:X\to Y$ be a quasi-perfect morphism between quasi-compact quasi-separated schemes. Then $$\Alg(\Lbf f^*):\Alg(\Dcal_{\QCoh}(Y))\to \Alg(\Dcal_{\QCoh}(X))$$ admits a left adjoint. If moreover $f$ is flat, then 
$$\Alg(f^*):\Alg(\QCoh(Y))\to \Alg(\QCoh(X))$$ admits a left adjoint.
\begin{proof}
Recall that a morphism of schemes $f$ is called \textit{quasi-perfect} if $\R f_*$ maps perfect complexes to perfect complexes.
We want to apply Corollary \ref{cor_adjoint} with $\Ccal:=\Dcal_{\QCoh}(Y)$ and $\Dcal:=\Dcal_{\QCoh}(X)$, and $F:=\Lbf f^*$. The $\infty$-categories $\Ccal,\Dcal$ are closed monoidal compactly generated and $F$ is monoidal and admits a right adjoint by $\R f_*$. As $f$ is itself quasi-compact quasi-separated, the functor $\R f_*$ maps $\Dcal_{\QCoh}(X)$ to $\Dcal_{\QCoh}(X)$). References for other conditions of Corollary \ref{cor_adjoint}:
\begin{enumerate}
\setcounter{enumi}{1}
    \item By \cite[\href{https://stacks.math.columbia.edu/tag/09M1}{Tag 09M1}]{stack}, the full $\infty$-subcategory of compact objects of quasi-coherent derived category coincides with the $\infty$-category of perfect complexes. Thus $\Ccal^{\omega}$ contains $\Ocal_Y$, is closed under derived tensor product (\cite[\href{https://stacks.math.columbia.edu/tag/09J5}{Tag 09J5}]{stack}), and is closed under dual by \cite[\href{https://stacks.math.columbia.edu/tag/08DQ}{Tag 08QD}]{stack}, which also says the natural transformation $\id_{\Ccal^{\omega}} \Rightarrow \D_{\Ccal^{\omega}}^{\op}\D_{\Ccal^{\omega}}$ are isomorphic. We have similar results for $\Dcal;$
    \item See \cite[\href{https://stacks.math.columbia.edu/tag/09UA}{Tag 09UA}]{stack} and recall the definition of quasi-perfect morphism;
    \item See \cite[\href{https://stacks.math.columbia.edu/tag/0GM7}{Tag 0GM7}]{stack}.
\end{enumerate}
If $f$ is flat then $f^*=\Lbf f^*$ so $\Lbf f^*$ is both left and right t-exact. Then Lemma \ref{lr_t-exact} applies that $\Alg(f^*)$ admits a left adjoint.
\end{proof}
\end{theorem}
\begin{remark}\label{rmk:module_lef}
In the situation of Theorem \ref{ad_module}, we actually have a left adjoint of $f^*:\QCoh(Y)\to \QCoh(X)$ which we will denote by $f_+$ in next section. Indeed, let us be in the situation of Lemma \ref{lr_t-exact}. Similar to the proof of Lemma \ref{lr_t-exact}, we have $\pi_0^{\Ccal}H|_{\Dcal_{\geq 0}}i_\Dcal$ is the left adjoint of $F|_{\Ccal^\heart}$ as for $M\in \Dcal^\heart$ and $N\in \Ccal^\heart$, one has 
$$
\begin{matrix*}
\Hom_{\Ccal^\heart}(\pi_0^{\Ccal}H|_{\Dcal_{\geq 0}}i_\Dcal M,N)&\simeq&\Hom_{\Ccal_{\geq 0}}(H|_{\Dcal_{\geq 0}}i_\Dcal M,i_\Ccal N)\\
&\simeq &\Hom_{\Dcal_{\geq 0}}(i_\Dcal M,F|_{\Ccal_{\geq 0}}i_\Ccal N)\\
&\simeq & \Hom_{\Dcal_{\geq 0}}(i_\Dcal M,i_\Dcal F|_{\Ccal^\heart} N)\\
&\simeq& \Hom_{\Dcal^\heart}(\pi_0^{\Dcal}i_\Dcal M,F|_{\Ccal^\heart} N)\\
&\simeq& \Hom_{\Dcal^\heart}(M,F|_{\Ccal^\heart} N).
\end{matrix*}
$$
Applying this formula of left adjoint to $F|_{\Ccal^\heart}=f^*,$ we obtain a formula $$f_+(\Fcal)\simeq \Hrm_0(\R f_*(\Fcal^\vee))^\vee$$ where $\Fcal$ is a perfect complex and $(-)^\vee = \R \Hom_{\Ocal_Y}(-,\Ocal_Y).$ 
\end{remark}
\section{A category of weighted polynomial laws}
\subsection{Polynomial laws over schemes}
We start by collecting some global definitions and results on polynomial laws. All of the proofs are direct consequences of \cite{Roby_1963}. 

For convenience in developing our theory as for some reasons such as the usual pullback functor of quasi-coherent sheaves only preserve linear morphisms, we adopt the language of functors of points which is developed in \cite[\textsection 1, $\text{n}^{\text{o}}$ 4]{Demazure_1970}. Here, we follow \cite[Part 2, Section 8]{Karpenko_2016}. Considering a scheme $S\to \Spec(R)$ (we only consider $R=\Zbb$) as a functor $\CAlg(R)\to \Set$ given by $X(A):=\Hom_{\Spec(R)}(\Spec(A),X)$ for any commutative $R$-algebra $A$.  The important convention is that a presheaf of $\Ocal_S$-modules is a functor $\CAlg(S)\to \Set$ which sends each object in $\CAlg(S)$ to a set equipped with a module structure. For reader's convenience, we recall their definitions:
\begin{definition}
Let $S:\CAlg(R)\to \Set$ be a scheme over $\Spec(R)$. The category $\CAlg(S)$ is defined as follows:
\begin{enumerate}
    \item The objects are the pairs $(A,x)$ where $A\in \CAlg(R)$ and $x\in X(A)$;
    \item A morphism between $(A,x)$ and $(B,y)$ is a homomorphism $f:A\to B$ of $R$-algebras such that $X(f)(x)=y.$
\end{enumerate}
An $S$-functor is then a functor $\CAlg(S)\to \Set$.
\end{definition}

\begin{definition}\label{def:qch}
\begin{enumerate}
    \item An $\Fcal$-functor $\Gcal$ is called a \textit{presheaf of modules} if for any $(A,x)\in \CAlg(\Fcal)$, we can put an $A$-module structure on the set $\Gcal(A,x)$ so that for any morphism $f:(A,x)\to (B,y),$ the morphism $\Gcal(A,x)\to \Gcal(B,y)$ is a homomorphism of $A$-modules. 
    \item Let $\Gcal$ be a presheaf of modules. We say that $\Gcal$ is a \textit{quasi-coherent sheaf} if the morphisms $\Gcal(A,x)\to \Gcal(B,y)$ induce isomorphisms $\Gcal(A,x)\otimes_A B\to \Gcal(B,y)$ of $B$-modules.
\end{enumerate}
\end{definition}
\begin{remark}
A nice feature of Definition \ref{def:qch}(2) is that we do not require the sheaf condition. This saves us some space in the following discussion.
\end{remark}
We now give some definitions and results parallel to \cite{Roby_1963}.
\begin{definition}\label{def:law_sh}
Let $\Fcal$ and $\Gcal$ be quasi-coherent sheaves over a scheme $S$. 
\begin{enumerate}
    \item A \textit{polynomial law} $\psi$ from $\Fcal$ to $\Gcal$ is a natural transformation $\Fcal\Rightarrow \Gcal$ as presheaves of sets. We will usually write $\psi:\Fcal \leadsto \Gcal$. By definition, a polynomial law $\psi:\Fcal \leadsto \Gcal$ is a datum consisting of maps of sets $$\{\psi_{(A,x)}:\Fcal(A,x)\to\Gcal(A,x)\}_{(A,x:\Spec(A)\to X)\in \CAlg(S)},$$ such that for every morphism $f:(A,x)\to (B,y)$, the following diagram
$$
\xymatrix{
\Fcal(A,x) \ar[r]^{\psi_{(A,x)}} \ar[d]_{\Fcal(f)}& 
\Gcal(A,x) \ar[d]^{\Gcal(f)}\\
\Fcal(B,y) \ar[r]_{\psi_{(B,y)}} & \Gcal(B,y)
}
$$
commutes. We let $\Prm_S(\Fcal,\Gcal)$ denote the set of all polynomial laws.
\item Let $\psi:\Fcal\leadsto\Gcal$ be a polynomial law of quasi-coherent sheaves over a scheme $S$. We say that $\psi$ is \textit{homogeneous of degree $d$} if for every $z\in \Fcal(A,x)$ and $r\in A$, we have 
$$\psi_{(A,x)}(zr)=\psi_{(A,x)}(z)r^d.$$ 
We let $\Prm^d_S(\Fcal,\Gcal)$ denote the set of all homogeneous polynomial laws of degree $d$. 
\end{enumerate}
\begin{remark}
Let $(A,x)\in \CAlg(S)$, then $\psi$ restricts to a polynomial law $$\psi|_{(A,x)}:\Fcal(A,x)\leadsto \Gcal(A,x)$$ as follows. For any algebra $f:A\to B$, let $y=(X(A)\to X(B))(x)$ and define $(\psi|_{(A,x)})_B$ as the dashed map making the following diagram commute:
$$
\xymatrix{
\Fcal(A,x)\otimes B \ar@{-->}[rr]^{} \ar[d]&& \Gcal(A,x) \otimes B \ar[d]\\
\Fcal(B,y) \ar[rr]_{\psi_{(B,y)}} && \Gcal(B,y),
}
$$
which is possible since the vertical maps are isomorphisms.
\end{remark}
The following characterization of polynomial laws is similar to \cite[Theorem I.1]{Roby_1963}:
\begin{theorem}
If $\psi:\Fcal\leadsto \Gcal$ is a polynomial law on scheme $S$ then for any integer $n\geq 0$, there exists a uniquely morphism $\alpha:\Ocal_X^{\oplus \Zbb^{\oplus n}}\to \Fcal$ satisfying the following conditions:
\begin{enumerate}
    \item $\alpha$ is locally finite, i.e., for any $(A,x)\in \CAlg(S),$ the family $(\alpha_{(A,x)}(e_i))_{i\in \Zbb^{\oplus n}}$ has only a finite number of non-zero elements, where $(e_i)_{i\in \Zbb^{\oplus n}}$ is the canonical basis of $A^{\oplus \Zbb^{\oplus n}}$.
    \item For any $(A,x)\in \CAlg(S)$ and any $n$ elements $m_1,\dots,m_n\in \Fcal(A,x),$ one has:
$$\psi_{(A,x)}(m_1\otimes r_1+\dots+m_n\otimes r_n)=\sum_{(k_1,\dots,k_n)\in \Zbb^{\oplus n}}\alpha_{(A,x)}(e_{(k_1,\dots,k_n)})\otimes r_1^{k_1}\cdots r_n^{k_n}$$
for any $r_1,\dots,r_n\in A.$
\end{enumerate}
\end{theorem}
\end{definition}
In a similar manner to \cite[Section I.5]{Roby_1963}, we define sum of polynomial laws as follows:
\begin{definition}
Let $(\psi_i:\Fcal\leadsto \Gcal)_{i\geq 0}$ be a locally finite family of  polynomial laws over $S$, i.e., for any $(A,x)\in \CAlg(S)$ and any $z\in \Fcal(A,x),$ the family $(\psi_{i,(A,x)}(z))_{i\geq 0}$ has only a finite number of of non-zero elements. Then we define a polynomial law $\sum_i \psi_i:\Fcal\leadsto \Gcal$ by the formula
$$\left(\sum_i \psi_i\right)_{(A,x)}(z)=\sum_i \psi_{i,(A,x)}(z)$$
for $(A,x)\in \CAlg(S)$ and $z\in \Fcal(A,x).$
\end{definition}
\begin{proposition}\label{prop:decomposition}
Let $\psi:\Fcal\leadsto \Gcal$ be a polynomial law over $S$. Then there there are unique homogeneous component laws $\psi_i$ of degree $i$ such that $\psi=\sum_{i\geq 0}\psi_i.$
\begin{proof}
For each $(A,x)\in \CAlg(S),$ we have homogeneous component polynomial laws $(\psi_{i,(A,x)})$ of $\psi|_{(A,x)}$ by \cite[Proposition I.4]{Roby_1963}. Then $\psi_i$ is the unique homogeneous polynomial laws which satifies $\psi_i|_{(A,x)}=\psi_{i,(A,x)}.$
\end{proof}
\end{proposition}
\begin{proposition}\label{prop:corr_sheaf}
Let $\Fcal$ be a quasi-coherent sheaf over a scheme $S$. Then there is a polynomial law $\lrm^d_\Gcal:\Gcal\leadsto \Gamma^d(\Gcal)$ such that the precomposition induces a bijection between $\Hom_{\Ocal_S}(\Gamma^d(\Fcal),\Gcal)$ and $\Prm^d_S(\Fcal,\Gcal)$.  
\end{proposition}
\begin{proof}
Recall that by {\cite[Proposition IV.1]{Roby_1963}\label{prop:correspond}}, for any $R$-module $M,$ the maps 
$$\lrm^{d}_{M,A}:M\otimes A \to \Gamma^{d}(M)\otimes A,\quad  x\otimes r \mapsto x^{[d]}\otimes r^d$$ define a homogeneous polynomial law of degree $d$ on $(M,\Gamma^d(M))$ and there exists a canonical isomorphism 
$$\Hom_R(\Gamma^d(M),N) \to \Prm^d(M,N), \quad \phi \mapsto \{(\phi\otimes A)\circ \lrm^{d}_{M,A}\}_{A\in \CAlg(R)}.$$
For each $(A,x)\in \CAlg(S),$ we define $$\lrm^{d}_{\Gcal,(A,x)}:\Gcal(A,x)\to\Gamma^d(\Gcal(A,x))$$ by $z\mapsto \lrm^d_{\Gcal(A,x)}(z)$. This gives the desired polynomial law $\lrm^d_{S}: \Gcal \leadsto \Gamma^d(\Gcal)$. 
\end{proof}
\begin{proposition}\label{prop:can_gamma}
There is a canonical injection $\Prm_X(\Fcal,\Gcal)\hookrightarrow\Hom_{\Ocal_X}(\Gamma(\Fcal),\Gcal).$ 
\begin{proof}
The map is defined as follows. Let $\psi:\Fcal\leadsto \Gcal$ be any polynomial law. By Proposition \ref{prop:decomposition}, we have a decomposition $\psi=\sum_{i\geq 0}\psi_i$. By Proposition \ref{prop:corr_sheaf}, each $\psi$ induces a linear morphism $\Gamma^i(\Fcal)\to \Gcal$ and hence a linear morphism $\Gamma(\Fcal)\to \Gcal.$ It is obvious from Proposition \ref{prop:decomposition} and Proposition \ref{prop:corr_sheaf} this map is injective.
\end{proof}
\end{proposition}
We say that a polynomial law $\psi$ is \textit{of finite degree} if $\mathrm{argmax}_{i\geq 0}\psi_i<\infty.$ Now Proposition \ref{prop:can_gamma} induces the following isomorphism:
\begin{lemma}
Let $\Prm_S^{\mathrm{fin}}(\Fcal,\Gcal)$ be the set of polynomial laws of finite degree. Then
$$\Prm_S^{\mathrm{fin}}(\Fcal,\Gcal)\simeq \bigoplus_{i\geq 0} \Hom_{\Ocal_S}(\Gamma^{i}(\Fcal),\Gcal).$$
\end{lemma}
We introduce an immediate category that will be useful:
\begin{definition}
We define\textit{ the category of quasi-coherent sheaves with polynomial morphisms} $\Prm_S$ as follows:
\begin{enumerate}
    \item The objects are quasi-coherent sheaves over $S.$
    \item The set of morphisms between sheaves $\Fcal$ and $\Gcal$ is $\Prm_S(\Fcal,\Gcal)$ in Definition \ref{def:law_sh}.
\end{enumerate}
\end{definition}
We now want to define a pullback functor and a pushforward functor between these categories. Note that the pullback functor of modules do not do not give rise such a functor as it only acts on linear momorphisms so we prefer the following direct definition:
\begin{definition}
Let $f:Y\to S$ be a morphism of schemes and let $\psi:\Fcal\leadsto \Gcal$ be a polynomial law over $S$. Then the pullback \textit{pullback polynomial law} $f^*\psi$ is given by the formula
$$f^*\psi(A,x)=\psi(A,f_A(x)),$$
for any $(A,x)\in \CAlg_Y.$ It is obvious that $f^*(\psi\circ \phi) = f^*(\psi)\circ f^*({\phi})$, so this extends to a functor $f^*:\Prm_S\to \Prm_Y.$   
\end{definition}
\begin{definition}
Let $f:Y\to S$ be quasi-compact quasi-separated morphism and $\psi:\Fcal\leadsto \Gcal$ be a polynomial law over $Y$. Then the pushforward functor define the \textit{pushforward polynomial law} $f_*\psi$ since $f_*$ maps quasi-coherent sheaves to quasi-coherent sheaves (\cite[\href{https://stacks.math.columbia.edu/tag/01LC}{Tag 01LC}]{stack}). This extends to a functor  $f^*:\Prm_Y\to \Prm_S.$   
\end{definition}
\subsection{Weighted polynomia laws}
\begin{definition}\label{def:w_poly}
Let $d$ be a non-negative integer. The category $\poly^d(S)$ is defined as follows:
\begin{itemize}
    \item The objects are tuples $\Psi=(\Fcal, \Gcal, \Acal_0, \Acal_1,\dots, \Acal_d, \psi_{-1}, \psi_0, \psi_1, \dots, \psi_d)$ where $\Fcal,\Gcal\in \QCoh_S,\ \Acal_i\in \Alg_S$ for $i=0, 1,\dots,d$;
     $\psi_{-1}:\Fcal\leadsto \Gcal$ is a polynomial law, and each $\psi_i:\Fcal\leadsto \Acal_i$ is a polynomial law for $i=0,1,\dots,d$.
     \item A morphism from $\Psi$ to $\Psi'$ is a tuple $(u,u_{-2},u_{-1},u_0,\dots,u_d)$ where 
    \begin{itemize}
        \item $u:\Fcal\leadsto \Fcal'$ is a polynomial law, and $u_{-1}: \Gcal\to \Gcal'$ is a linear morphism;
        \item Each $u_i:\Acal_i\to \Acal'_i$ is a morphism of sheaves of algebras for $i=-0,1,\dots,d$;
\end{itemize}
such that the following diagrams 
$$
\xymatrix{
\Fcal\ar@{~>}[r]^{\psi_{-1}}\ar@{~>}[d]_{u} & \Gcal\ar[d]^{u_{-1}}\\
\Fcal'\ar@{~>}[r]_{\psi'_{-1}} & \Gcal'
}
\xymatrix{
\Fcal \ar@{~>}[r]^{\psi_i} \ar@{~>}[d]_{u}& \Acal_i\ar[d]^{u_i}\\
\Fcal' \ar@{~>}[r]_{\psi'_i} & \Acal'_i
}
$$
commute for $i=0,1,\dots,d$. 
\end{itemize}
If $\psi_{-1},\psi_0,\dots,\psi_{d-1}$ are of finite degree then we say that $\Psi$ is \textit{of finite degree}. We let $\poly^d_{\fin}(S)$ denote the subcategory obtained by restricting the collection of morphisms to those tuples $(u,(u_i)_{-1\leq i\leq d})$ where $u$ is of finite degree. When $S=\Spec(R),$ we write $\poly^d(R), \poly^d_{\fin}(R)$ for $\poly^d(S),\poly^d_{\fin}(S),$ respectively.

We also consider the following more general variant of the above category which includes bilinear forms. The category $\poly^d_{\mathrm{bi}}(S)$ is defined as follows:
\begin{itemize}
    \item The objects are tuples $\Psi=(\Fcal, \Acal_{-2}, \Gcal, \Acal_0, \Acal_1,\dots, \Acal_d, \psi_{-2}, \psi_{-1}, \psi_0, \psi_1, \dots, \psi_d)$ where $\Fcal,\Gcal\in \QCoh_S,\ \Acal_i\in \Alg_S$ for $i=-2, 0,\dots,d$; $\psi_{-2}:\Gcal\otimes \Gcal\to \Acal_{-2}$ is a linear morphism,
     $\psi_{-1}:\Fcal\leadsto \Gcal$ is a polynomial law, and each $\psi_i:\Fcal\leadsto \Acal_i$ is a polynomial law for $i=0,1,\dots,d$.    
    \item A morphism from $\Psi$ to $\Psi'$ is a tuple $(u,u_{-2},u_{-1},u_0,\dots,u_d)$ where 
    \begin{itemize}
        \item $u:\Fcal\leadsto \Fcal'$ is a polynomial law, and $u_{-1}: \Gcal\to \Gcal'$ is a linear morphism;
        \item Each $u_i:\Acal_i\to \Acal'_i$ is a morphism of sheaves of algebras for $i=-2,0,\dots,d$;
    \end{itemize}
such that the following diagrams 
$$
\xymatrix{
\Gcal\otimes \Gcal \ar[r]^{\psi_{-2}}\ar[d]_{u_{-1}\otimes u_{-1}} & \Acal_{-2}\ar[d]^{u_{-2}}\\
\Gcal'\otimes \Gcal'\ar[r]_{\psi'_{-2}} & \Acal'_{-2}
}
\xymatrix{
\Fcal\ar@{~>}[r]^{\psi_{-1}}\ar@{~>}[d]_{u} & \Gcal\ar[d]^{u_{-1}}\\
\Fcal'\ar@{~>}[r]_{\psi'_{-1}} & \Gcal'
}
\xymatrix{
\Fcal \ar@{~>}[r]^{\psi_i} \ar@{~>}[d]_{u}& \Acal_i\ar[d]^{u_i}\\
\Fcal' \ar@{~>}[r]_{\psi'_i} & \Acal'_i
}
$$
commute for $i=0,1,\dots,d$. 
\end{itemize}
The categories $\poly^d_{\mathrm{bi},\fin}(S)$, $\poly^d_{\mathrm{bi}}(R), \poly^d_{\mathrm{bi},\fin}(R)$ are defined similarly. We refer to objects of the categories as\textit{ weighted polynomial laws} and its polynomial laws as \textit{weights}.
\end{definition}

We now consider a functor from the category of algebras to the category of weighted polynomial laws:
\begin{definition}
We let $F^d_S:\Alg_S\to \poly^d_{\mathrm{bi}}(S)$ be the functor that maps $(\Acal\in \Alg_S)$ to the weighted polynomial law $$(\Acal^{\oplus (d+1)},\Acal,\dots,\Acal,\pr_{-2},\pr_{-1}, \pr_{0},\dots,\pr_{d}),$$
where for any $(A,x)\in \CAlg(S)$, we define
\begin{itemize}
\item$\pr_{-2,(A,x)}:\Acal(A,x)\otimes \Acal(A,x) \to \Acal(A,x)$ sends $x\otimes y$ to $xy-yx.$
    \item Each $\pr_{i,(A,x)}:\Acal(A,x)^{\oplus(d+1)}\to \Acal(A,x)$ is the $(i+2)$-th coordinate projection for $i=-1,0,\dots,d-1.$
    \item $\pr_{d,(A,x)}:\Acal(A,x)^{\oplus(d+1)}\to\Acal(A,x)$ sends $(a_{-1},\dots,a_{d-1})$ to $a_0a_{-1}^d+a_1a_{-1}^{d-1}+\dots+a_{d-1}a_{-1}.$
\end{itemize}
\end{definition}
\begin{remark}\label{rmk:compatible}
The collection of functors $(F^d_S:\Alg_S\to \poly^d_{\mathrm{bi}}(S))_{S}$ is \textit{compatible with base change}: for any morphism $f:S'\to S$ with schemes $S,S'$, we have a natural isomorphism $f^*F_S\simeq F_{S'}f^*$ where $f^*$ on the left is pullback functor of polynomial law and $f^*$ on the right is pullback functor of sheaves of algebras.
\end{remark}
We now give the definition of the \textit{Clifford algebra functor} as the left adjoint of functor $F^d_S.$
\begin{theorem}\label{thm:Cliff}
The functor $F^d_S:\Alg_S\to \poly^d_{\mathrm{bi}}(S)$ admits a left adjoint. More explicitly, let denote this left adjoint by $\Cl$ and let $\Psi=(\Fcal,\Acal_{-2},\Gcal,\Acal_{0},\dots,\Acal_d,\psi_{-2},\dots,\psi_d)$ be a weighted polynomial law. Then $\Cl(\Psi)$ fits into the following  colimit diagram
$$
\xymatrix{
\T(\Gcal\otimes \Gcal)\ar[r]\ar[d] & \Acal_0\ast\cdots\ast\Acal_{d-1} * \T(\Gcal) \ar[d]&\T(\Gamma(\Fcal))\ar[l] \ar[d] \\
\Acal_{-2}\ar@{-->}[r]& \Cl(\Psi)& \Acal_d\ar@{-->}[l]
}
$$
where $\Acal_0\ast\cdots\ast\Acal_{d-1} * \T(\Gcal)$ is the coproduct of sheaves of algebras. The morphisms in the diagram is given as follows:
\begin{enumerate}
    \item The morphism $\T(\Gcal \otimes \Gcal)\to \Acal_0\ast\cdots\ast\Acal_{d-1} * \T(\Gcal)$ is induced by linear morphism $\Gcal\otimes \Gcal\to \T(\Gcal),m\otimes n \mapsto m\otimes n - n\otimes m.$
    \item The morphism $\T(\Gamma(\Fcal))\to \Acal_d$ is induced by the law $\psi_d$ and the morphism $\T(\Gamma(\Fcal))\to \Acal_0\ast\cdots\ast\Acal_{d-1} * \T(\Gcal)$ is induced by the law $ \psi_{-1}^d\psi_0+\psi_{-1}^{d-1}\psi_1+\cdots+\psi_{-1}\psi_{d-1}:
    \Fcal \leadsto \Acal_0\ast\cdots\ast\Acal_{d-1} * \T(\Gcal), $
\end{enumerate}
\end{theorem}
\begin{proof}
Noting that $\Alg_S$ has all small colimits so let $\Acal$ be the colimit of the diagram:
$$
\xymatrix{
\T(\Gcal\otimes \Gcal)\ar[r]\ar[d] & \Acal_0\ast\cdots\ast\Acal_{d-1} * \T(\Gcal) \ar[d]&\T(\Gamma(\Fcal))\ar[l] \ar[d] \\
\Acal_{-2}\ar@{-->}[r]& \Acal& \Acal_d\ar@{-->}[l]
}
$$
Then $\Hom_{\Alg_S}(\Acal,\Acal')$ coincides with the set of tuples $(u_{-2},\dots,u_d)$ where 
\begin{enumerate}
    \item $u_{-1}:\Gcal\to\Acal$ is a linear morphism (induced by the morphism $\Acal_0\ast\cdots\ast\Acal_{d-1}*\T(\Gcal)$ in $\Alg_S$);
    \item Each $u_{i}:\Acal_{i}\to\Acal$ is a morphisms of algebras for $i=-2,0,\dots,d$;
\end{enumerate}
Which make the following diagrams commute:
 $$\xymatrix{
\Gcal \otimes \Gcal \ar[rr]^-{\psi_{-2}} \ar[d]_{u_{-1}\otimes u_{-1}}& &\Acal_{-2} \ar[d]^{u_{-2}}\\
\Acal \otimes \Acal \ar[rr]_-{x\otimes y \mapsto xy-yx} & &\Acal
    },
    \xymatrix{
    \Acal_i \ar[d]_{u_i}& \Fcal \ar[l]_{\psi_i}\ar[d]\ar[r]^{\psi_{d}}& \Acal_d\ar[d]^{u_{d}}\\
    \Acal&\Acal^{\oplus d+1}  \ar[l]^-{\pr_i}\ar[r]_-{\pr_{d}} &\Acal
    }$$
    because of the commutativity of the diagrams: $$\xymatrix{
    \T(\Gcal\otimes \Gcal)\ar[d] \ar[r] & \Acal_0\ast \cdots\ast \Acal_{d-1}\ast \T(\Gcal) \ar[d]\\
    \Acal_{-2} \ar[r] & \Acal},
    \xymatrix{
    \T(\Gamma(\Fcal)) \ar[d] \ar[r]& \Acal_d\ar@{-->}[d]\\
    \Acal_0\ast\cdots\ast\Acal_{d-1} * \T(\Gcal)\ar@{-->}[r]& \Acal,
    }$$
where we use Proposition \ref{prop:can_gamma} for the second diagram. Hence these tuples correspond to the morphisms from $\Psi$ to $F(\Acal)$ in $\poly^d_{\mathrm{bi}}(S)$ (note that the law $\Fcal\leadsto \Acal^{\oplus d+1}$ is completely determined by $u_i$ as it is given by $(u_{-1}\circ \psi_{-1},\cdots,u_{d-1}\circ \psi_{d-1})$.
\end{proof}
\begin{remark}
In the explicit construction of $\Cl(\Psi)$ in Theorem \ref{thm:Cliff}, we can replace $\Gamma(\Fcal)$ by a finite direct sum of $\Gamma^j(\Fcal)$ where $j$ is among degrees of homogeneous component laws of $\psi_{-1}^d\psi_0+\psi_{-1}^{d-1}\psi_1+\cdots+\psi_{-1}\psi_{d-1}$ and $\psi_d$ if they are of finite degree as in following examples.  
\end{remark}
\begin{example}
Any linear morphism $\psi_{-2}=\psi:M\otimes M\to R$ determines a weighted polynomial law $\Psi$ in $\poly^1(R)$ given by the following weights:
\begin{itemize}
    \item $\psi_{-1}:M\to M$ is the identity.
    \item $\psi_0: M\to R$ is the zero map.
    \item $\psi_1:M\to \T(M)$ is the canonical linear inclusion.
\end{itemize}
By Theorem \ref{thm:Cliff}, the algebra $\Cl(\Psi)$ fits into the following colimit diagram:
$$
\xymatrix{
\T(M\otimes M)\ar[r]\ar[d] & R\ast \T(M)\ar[d]&\T(M)\ar[l] \ar[d] \\
R\ar@{-->}[r]& \Cl(\Psi)& \T(M).\ar@{-->}[l]
}
$$
Note that $R*\T(M)\simeq \T(M)$ and the map $\T(M)\to R*\T(M)$ maps $M$ to $\{0\}$, which implies that the map $\T(M)\to \Cl(\Psi)$always  maps $M$ to $\{0\}$. Hence $\Cl(\Psi)$ is the pushout of the following diagram:
$$
\xymatrix{
\T(M\otimes M) \ar[r]\ar[d] &\T(M) \ar@{-->}[d]\\
R \ar@{-->}[r] & \Cl(\Psi).
}
$$
Therefore $\Cl(\Psi)\simeq \T(M)/(m\otimes n-n\otimes m-\psi(m,n)).$ The Weyl algebra is an instance of this algebra when $M$ is a vector space of even dimension and $\psi$ is a symplectic form. 
\end{example}
\begin{example}
Let $q:M\to R$ be a quadratic form. Then $q$ induces a homogeneous polynomial law $\psi_2:M\leadsto R$ of degree $2$. This gives a weighted polynomial law $\Psi$ in $\poly^2(R)$ by picking the following weights:
\begin{itemize}
\item $\psi_{-2}:M\otimes M\to \T(M)$ is given by $z\otimes z'\mapsto z\otimes z'-z'\otimes z$.
    \item $\psi_{-1}:M\to M$ is the identity map.
    \item $\psi_0:M\to R$ maps everything to $\{1\}$.
    \item $\psi_1: M\to R$ is the zero map.
\end{itemize}
The by Theorem \ref{thm:Cliff}, the algebra $\Cl(\Psi)$ fits into the following colimit diagram:
$$
\xymatrix{
\T(M\otimes M)\ar[r]\ar[d] &R\ast R\ast \T(M)\ar@{-->}[d]&\T(\Gamma^2(M))\ar[l] \ar[d] \\
\T(M)\ar@{-->}[r]& \Cl(\Psi)& R.\ar@{-->}[l]
}
$$
Note that the left hand square always commutes for any morphism $\T(M)$ to $\Cl(\Psi)$ so $\Cl(\Psi)$ is the pushout of the square
$$
\xymatrix{
\T(\Gamma^2(M))\ar[d]\ar[r]& \T(M)\ar@{-->}[d]\\
R \ar@{-->}[r]& \Cl(\psi).
}
$$
Also by Theorem \ref{thm:Cliff}, we have the following universal property of $\Cl(\Psi):$ For any $R$-algebra $A,$ we have $\Hom_{\Alg_R}(\Cl(\Psi),A)$ is bijective to the set of tuples $(u_{-2},u_{-1},\dots,u_2)$ making the following diagrams commute:
$$
\xymatrix{
M\otimes M \ar[rrr]^{h\otimes k\mapsto h\otimes k-k\otimes h}\ar[d]_{u_{-1}\otimes u_{-1}} & &&\T(M)\ar[d]^{u_{-2}}\\
A\otimes A \ar[rrr]_{h\otimes k\mapsto hk-kh} && &A
}
\xymatrix{
M\ar[rr]^{\id_M}\ar@{~>}[d]_{u} & &M\ar[d]^{u_{-1}}\\
A^{\oplus 3}\ar[rr]_{(h,k,l)\mapsto h} && A
}
\xymatrix{
M \ar@{~>}[rr]^{1} \ar@{~>}[d]_{u}&& R\ar[d]^{u_0}\\
A^{\oplus 3} \ar[rr]_{(h,k,l)\mapsto k} && A
}
$$
$$
\xymatrix{
M \ar[rr]^{0} \ar@{~>}[d]_{u}&& R\ar[d]^{u_1}\\
A^{\oplus 3} \ar[rr]_{(h,k,l)\mapsto l} && A
}
\xymatrix{
M \ar@{~>}[rr]^{q} \ar@{~>}[d]_{u}&& R\ar[d]^{u_2}\\
A^{\oplus 3} \ar[rr]_{(h,k,l)\mapsto kh^2+lh} && A.
}
$$
Note that $u_0=u_1=u_2=i_A$ the structure map of the algebra $A$. The first diagram always commutes for any linear map $u_{-1}$. The polynomial law $u$ decomposes into polynomial law $v_1,v_2,v_3:M\leadsto A$. The second commutative diagram implies that $u_{-1}=v_1.$ The third one implies that $v_2=1.$ The fourth one implies that $v_3=0.$ Hence the last one implies that $\Hom_{\Alg_R}(\Cl(\Psi),A)$ is bijective to the set of linear morphism $u_{-1}\to A$ such that $i_A\circ q = v_2v_{1}^2+v_3v_{1}=u_{-1}^2.$ This is the universal property of quadratic Clifford algebra.
\end{example}
\begin{example}\label{ex:KL}
The weighted non-diagonal Clifford algebra of homogeneous polynomials \cite[Appendix A]{KL_2015} is defined as follows. Let $f_m,f_{2m},\dots,f_{dm}\in R[x_1,\dots,x_n]$ where $f_i$ is homogeneous of degree $i.$ Then the Clifford algebra $\Cl(f_m,\dots,f_{dm})$ is given by 
$$\Cl(f_m,\dots,f_{dm})=k\langle a_J\rangle_{|J|=m}/I$$
where $I$ is the ideal generated by the coefficients of the variables $x_i$ in the expression 
$$\left(\sum_{|J|=m} a_Jx^J\right)^d-\sum_{l=1}^d\left(\left(\sum_{|J|=m}a_Jx^J\right)^{d-l}f_{lm}(x_1,\dots,x_n)\right)\in R\langle a_J \rangle_{|J|=m}[x_1,\dots,x_n].$$
On the other hand, the polynomials $f_m,f_{2m},\dots,f_{dm}$ give us the a weighted polynomial law $\Psi$ in $\poly^d_{\mathrm{bi}}(R)$ defined by the following weights:
\begin{itemize}
    \item  $\psi_{-2}:M\otimes M\to \T(M)$ is given by $z\otimes z'\mapsto z\otimes z'-z'\otimes z$ where $M=R^n$.
    \item $\psi_{-1}:M\leadsto N, x_1e_1+\cdots+x_ne_n\mapsto \sum_{j_1+\cdots+j_n=m}x_1^{j_1}\cdots x_n^{j_n}e_1^{j_1}\cdots e_n^{j_n},$ where $N=\Sym^m(M)$ and $(e_i)_{1\leq i\leq n}$ is a basis of $M$.
    \item $\psi_0:M\to R$ is the polynomial law $1.$
    \item $\psi_i:M\leadsto R$ is induced by $-f_{im}$ for $i=1,2,\dots,d-1$ and $f_{dm}$ for $i=d.$
\end{itemize}
Then $\Cl(\Psi)$ is the colimit of the diagram
$$
\xymatrix{
\T(M\otimes M)\ar[r]\ar[d] & \T(\Sym^m(M)) \ar[d]&\T(\Gamma^{md}(M))\ar[l] \ar[d] \\
\T(M)\ar@{-->}[r]& \Cl(\Psi)& R.\ar@{-->}[l]
}
$$
by Theorem \ref{thm:Cliff}. As in previous examples, $\Cl(\Psi)$ is solely the pushout of the following square:
$$
\xymatrix{
\T(\Gamma^{md}(M)) \ar[r] \ar[d]& \T(\Sym^m(M))\ar@{-->}[d]\\
R \ar@{-->}[r] & \Cl(\Psi).
}
$$
Hence is isomorphic to the quotient
$$
\frac{R\langle e_1^{j_1}\cdots e_n^{j_n} \mid j_1+\cdots j_n = m\rangle}{I}
$$
where $I$ is generated by the image of the morphism $\widetilde{\psi}_{d}-(\widetilde{\psi_{-1}^d+\psi_{-1}^{d-1}\psi_1+\cdots+\psi_{-1}\psi_{d-1}})$
where $\widetilde{\psi}_{d}$ and $\widetilde{\psi_{-1}^d+\psi_{-1}^{d-1}\psi_1+\cdots+\psi_{-1}\psi_{d-1}}$ are corresponding linear morphisms $\Gamma^{md}(M)\to \Sym^m(M).$ In fact, if suffices to pick $I$ as $$\left(\widetilde{\psi}_{d}(z^{[md]})-\widetilde{\psi_{-1}^d+\psi_{-1}^{d-1}\psi_1+\cdots+\psi_{-1}\psi_{d-1}}(z^{[md]})\mid z\in M\right).$$
Indeed, let $I'$ be the above ideal and let $$A=\frac{R\langle e_1^{j_1}\cdots e_n^{j_n} \mid j_1+\cdots j_n = m\rangle}{I'}.$$ Then the compositions $(\T(\Sym^m(M))\rightarrow A)\circ \widetilde{\psi}_{d}\circ\lrm^{md}_{M}$ and  $$(\T(\Sym^m(M))\rightarrow A)\circ (\widetilde{\psi_{-1}^d+\psi_{-1}^{d-1}\psi_1+\cdots+\psi_{-1}\psi_{d-1}})\circ\lrm^{md}_{M}$$
are equal. By Proposition \ref{prop:corr_sheaf}, $$(\T(\Sym^m(M))\rightarrow A)\circ \widetilde{\psi}_{d}=(\T(\Sym^m(M))\rightarrow A)\circ (\widetilde{\psi_{-1}^d+\psi_{-1}^{d-1}\psi_1+\cdots+\psi_{-1}\psi_{d-1}}).$$
It follows that there is a morphism $A\to \Cl(\Psi)$ which is an isomorphism with the inverse the projection $A\to \frac{R\langle e_1^{j_1}\cdots e_n^{j_n} \mid j_1+\cdots j_n = m\rangle}{I}.$ Thus $\Cl(\Psi)$ is isomorphic to 
$$
\frac{R\langle e_1^{j_1}\cdots e_n^{j_n} \mid j_1+\cdots j_n = m\rangle}{I}
$$
where $I$ is generated by 
$$f_{dm}(x_1,\cdots,x_n)-\left(\sum_{|J|=m}x_1^{j_1}\cdots x_n^{j_n}e_1^{j_1}\cdots e_n^{j_n}\right)^d+\left(\sum_{|J|=m}x_1^{j_1}\cdots x_n^{j_n}e_1^{j_1}\cdots e_n^{j_n}\right)^{d-1}f_m(x_1,\dots,x_n)$$
$$+\cdots+\left(\sum_{|J|=m}x_1^{j_1}\cdots x_n^{j_n}e_1^{j_1}\cdots e_n^{j_n}\right)f_{(d-1)m}(x_1,\dots,x_n)$$
for $x_1,\dots,x_n\in R$. The correspondence $ e_1^{j_1}\cdots e_n^{j_n}\mapsto a_{(j_1,\dots,j_n)}$ gives a canonical morphism $\Cl(\Psi)\to \Cl(f_m,\dots,f_{dm}).$ It is an isomorphism, for example when $d=2, m=1, f_1=0$, but this is no true in general. For example, consider the case where $R=\Fbb_2, m=1, d=2, n=1,$ and $f_{2}=0,f_1=y.$ Then by looking at expression $(ax+by)^2-(ax+by)y$, we obtain $$\Cl(f_1,f_2)\simeq \frac{R\langle a,b\rangle}{(a^2,ab-a,ba,b^2-b)},$$ and $$\Cl(\Psi)\simeq \frac{R\langle a,b\rangle}{(a^2,a^2+b^2+ab+ba-a-b,b^2-b)}.$$ 
Observe that $\Ker(\Cl(\Psi)\to \Cl(f_1,f_2))\neq 0$ since  $ab-a \notin (a^2,ab-a,ba,b^2-b).$ Indeed, if it does then there are $p_1,p_2,p_3,p_4,p_5,p_6\in R\langle a,b\rangle$ such that 
$$ab-a=p_1a^2p_2+p_3(a^2+b^2+ab+ba-a-b)p_4+p_5(b^2-b)p_6.$$
Comparing the degrees of both sides, we obtain an equation where we can choose all $p_i$ such that they are of degree $0,$ i.e., all $p_i\in R.$ This leads to a contradiction. 
In the case $M$ is not free, we can use the following weights:
\begin{itemize}
    \item $\psi_{-2}:M\otimes M\to \T(M), z\otimes z' \mapsto z\otimes z'-z'\otimes z$.
    \item $\psi_{-1}:M\to \Gamma^m(M), z\mapsto z^{[m]}.$
    \item $\psi_i:M\to R$ is a homogeneous polynomial law of degree $im$ for $i=0,1,\dots,d.$
\end{itemize}
\end{example}

\begin{lemma}\label{lm:func_cliff}
Let $f:Y\to S$ be a morphism of schemes. Then $\Cl(f^*\Psi)\simeq f^*\Cl(\Psi)$ for any $\Psi\in \poly^d_{\bi}(S)$ or $\Psi\in \poly^d(S).$
\end{lemma}
\begin{proof}
We prove the statement of $\Psi\in \poly^d_{\bi}(S)$, the other is similar. By Theorem \ref{thm:Cliff}, the sheaf
$\Cl(\Psi)$ fits into the following  colimit diagram
$$
\xymatrix{
\T(\Gcal\otimes \Gcal)\ar[r]\ar[d] & \Acal_0\ast\cdots\ast\Acal_{d-1} * \T(\Gcal) \ar[d]&\T(\Gamma(\Fcal))\ar[l] \ar[d] \\
\Acal_{-2}\ar@{-->}[r]& \Cl(\Psi)& \Acal_d\ar@{-->}[l]
}
$$
Since $f^*$ is a left adjoint, it commutes with all colimits, and it also commutes with tensor powers and $\Gamma^i$. This gives the corresponding colimit diagram
$$
\xymatrix{
\T(f^*\Gcal\otimes f^*\Gcal)\ar[r]\ar[d] & f^*\Acal_0\ast\cdots\ast f^*\Acal_{d-1} * \T(f^*\Gcal) \ar@{-->}[d]&\T(\Gamma(f^*\Fcal))\ar[l] \ar[d] \\
f^*\Acal_{-2}\ar@{-->}[r]& f^*\Cl(\Psi)& f^*\Acal_d.\ar@{-->}[l]
}
$$
One can verify that the solid arrows are induced by $f^*\Psi,$ so this diagram in fact defines $\Cl(f^*\Psi).$
\end{proof}
\subsection{Lifting of weighted polynomial laws}
\begin{definition}
Let $f:Y \to S$ be a morphism of schemes and 
\begin{enumerate}
    \item Let $\Psi_S\in \poly^d_{\mathrm{bi}}(S), \Psi_Y\in \poly^d_{\mathrm{bi}}(Y)$. We say that \textit{$\Psi_Y$ is a lift of $\Psi_S$ along $f$} if $\Psi_S$ corepresents the functor $$\Hom_{\poly^{d}_{\bi, \fin}(Y)}(\Psi_Y,f^*(-)):\Alg_S\to \Set.$$ 
    \item  Let $\Psi_S\in \poly^d(S), \Psi_Y\in \poly^d(Y)$. We say that \textit{$\Psi_Y$ is a lift of $\Psi_S$ along $f$} if $\Psi_S$ corepresents the functor $$\Hom_{\poly^{d}_\fin(Y)}(\Psi_Y,f^*(-)):\Alg_S\to \Set.$$ 
\end{enumerate}
\end{definition}
We immediately have the following:
\begin{proposition}\label{prop:lift_alg}
    Let $f:Y\to S$ be a morphism of schemes such that $f^*:\Alg_S\to \Alg_Y$ admits a left adjoint $f_\sharp.$ Let $\Psi_S\in \poly^{d}(S)$ or $\Psi_S\in \poly^{d}_{\bi}(S)$. If $\Psi_Y$ is a lift of $\Psi_S$ along $f$ then $f_\sharp(\Cl(\Psi_Y))\simeq \Cl(\Psi_S)$.
\end{proposition}
\begin{proof}We consider the case $\Psi_S\in \poly^d(S),$ the case $\Psi_S\in \poly^d_{\bi}(S)$ is similar. We have
$$
\begin{matrix*}[l]
\Hom_{\Alg_S}(f_\sharp\Cl(\Psi_Y),\Acal) &\simeq& \Hom_{\Alg_Y}(\Cl(\Psi_Y),f^*(\Acal))\\
&\simeq& \Hom_{\poly^{d}(Y)}(\Psi_Y,F_Yf^*(\Acal))\\
&\simeq& \Hom_{\poly^{d}_\fin(Y)}(\Psi_Y,F_Yf^*(\Acal))\\
&\simeq& \Hom_{\poly^{d}_{\fin}(Y)}(\Psi_Y,f^*F_S(\Acal))\\
&\simeq& \Hom_{\poly^{d}_\fin(S)}(\Psi_S,F_S(\Acal))\\
&\simeq& \Hom_{\poly^{d}(S)}(\Psi_S,F_S(\Acal))\\
&\simeq& \Hom_{\Alg_S}(\Cl(\Psi_S),\Acal),
\end{matrix*}
$$    
where the third bijetion is due to Remark \ref{rmk:compatible}.
\end{proof}
We now have a result to produce lifts in ideal situation:
\begin{theorem}\label{thm:main_compute}
Assume that $f:Y\to S$ is a flat quasi-perfect morphism between quasi-compact quasi-separated schemes. Let $f_+$ be the left adjoint of $f^*$.  Consider the following conditions on $\Fcal, \Gcal \in \QCoh(Y).$
\begin{enumerate}
    \item The morphism induced from oplax monoidal structure $f_+(\Gcal\otimes \Gcal)\to f_+(\Gcal)\otimes f_+(\Gcal)$ is an isomorphism.
    \item The morphism $f_+\Gamma^i(\Fcal)\to \Gamma^i(f_+\Fcal)$ defined as the composition $$f_+\Gamma^i(\Fcal)\xrightarrow{f_+\Gamma^i(\eta_{\Fcal})} f_+\Gamma^i(f^*f_+\Fcal)\to f_+f^*\Gamma^i(f_+\Fcal)\xrightarrow{\epsilon_{\Gamma^i(f_+\Fcal)}} \Gamma^i(f_+\Fcal)$$ is an isomorphism for all $i$, where $\eta$ is the unit of the adjunction $(f_+\dashv f^*)$.
\end{enumerate} 
\begin{itemize}
    \item Let $$\Psi=(\Fcal, \Acal_{-2}, \Gcal, \Acal_0, \Acal_1,\dots, \Acal_d, \psi_{-2}, \psi_{-1}, \psi_0, \psi_1, \dots, \psi_d)$$ be a weighted polynomial law of finite degree in $\poly^d_{\mathrm{bi}}(Y)$. Assume that $\Psi$ satisfies the above conditions then $\Psi$ is a lift along $f$. 
    \item Let $$\Psi=(\Fcal, \Gcal, \Acal_0, \Acal_1,\dots, \Acal_d, \psi_{-1}, \psi_0, \psi_1, \dots, \psi_d)$$ be a weighted polynomial law of finite degree in $\poly^d(Y)$. Assume that $\Psi$ satisfies the above condition (2) then $\Psi$ is a lift along $f$. 
\end{itemize}
\begin{proof}
We only consider the case $\Psi\in \poly^d_{\bi}(Y),$ the case $\Psi\in \poly^d(Y)$ is similar. With notations as in the statement, we first show that 
$$\Prm_{Y,\fin}(\Fcal,f^*(-))\simeq \Prm_{S,\fin}(f_+\Fcal,-).$$
Indeed, 
$$
\begin{matrix*}[l]
\Prm_{Y,\fin}(\Fcal,f^*(-))&\simeq& \bigoplus_{i\geq 0}\Hom_{\Ocal_Y}(\Gamma^i(\Fcal),f^*(-))\\
&\simeq&\bigoplus \Hom_{\Ocal_S}(f_+\Gamma^i(\Fcal),-)\\
&\simeq&\bigoplus \Hom_{\Ocal_S}(\Gamma^i(f_+\Fcal),-)\\
&\simeq& \Prm_{S,\fin}(f_+\Fcal,-).
\end{matrix*}
$$
Let $f_\sharp$ be the left adjoint of $f^*$. For any $\Acal\in \Alg_S$, composition with unit $\Acal\to f^*f_\sharp \Acal$ gives a map $$\Prm_{Y,\fin}(\Fcal,\Acal)\to \Prm_{Y,\fin}(\Fcal,f^*f_\sharp \Acal)\simeq \Prm_{S,\fin}(f_+\Fcal,f_\sharp \Acal).$$
We also have 
$$
\begin{matrix*}[l]
\Hom_{\Ocal_Y}(\Gcal\otimes \Gcal,f^*(-))&\simeq&  \Hom_{\Ocal_S}(f_+(\Gcal\otimes \Gcal),-)\\
&\simeq& \Hom_{\Ocal_S}(f_+\Gcal\otimes f_+\Gcal,-),
\end{matrix*}
$$
which for any $\Acal\in \Alg_S$ induces a map by compsition with unit $\Acal\to f^*f_\sharp \Acal$
$$\Hom_{\Ocal_Y}(\Gcal\otimes \Gcal,\Acal)\to \Hom_{\Ocal_Y}(\Gcal\otimes \Gcal,f^*f_\sharp\Acal)\simeq\Hom_{\Ocal_Y}(f_+\Gcal\otimes f_+\Gcal,f_\sharp \Acal)$$
Thus any $\Psi$ satisfying the condition in the statement gives rise to a $\Psi'\in \poly^d_{\mathrm{bi}}(S)$ with $$\Psi'=(f_+\Fcal, f_\sharp\Acal_{-2}, f_+\Gcal, f_\sharp\Acal_0, f_\sharp\Acal_1,\dots, f_\sharp\Acal_d, \psi'_{-2}, \psi'_{-1}, \psi'_0, \psi'_1, \dots, \psi'_d)$$
and $\psi'_i$ is induced by $\psi$. We now show that $\Psi$ is a lift of $\Psi'.$ Indeed, let $$\Psi''=(\Fcal'', \Acal''_{-2}, \Gcal'', \Acal''_0, \Acal''_1,\dots, \Acal''_d, \psi''_{-2}, \psi''_{-1}, \psi''_0, \psi''_1, \dots, \psi''_d)$$ be a weighted polynomial law in $\poly^d_{\mathrm{bi}}(S),$ then $\Hom_{\poly^d_{\mathrm{bi}}(Y)}(\Psi,f^*(\Psi''))$ is the set of tuples $\{u,u_{-2},\dots,u_d\}$ making the following diagrams, for $i=0,\dots,d$
$$\xymatrix{
\Gcal \otimes \Gcal \ar[r]^-{\psi_{-2}} \ar[d]_{u_{-1}\otimes u_{-1}}& \Acal_{-2} \ar[d]^{u_{-2}}\\
f^*\Gcal'' \otimes f^*\Gcal'' \ar[r]_-{f^*\psi''_{-2}}  &f^*\Acal''_{-2}
    },
    \xymatrix{
     \Gcal\ar[d]_{u_{-1}} &\Fcal\ar[l]_{\psi_{-1}} \ar[d]^{u}\ar[r]^{\psi_{i}}& \Acal_i\ar[d]^{u_{i}}\\f^*\Gcal''&f^*\Fcal''\ar[l]_{f^*\psi''_{-1}}\ar[r]^-{f^*\psi''_{i}} &f^*\Acal''_i
    }$$
commute. We now analyze each commutative diagram:
\begin{itemize}
    \item We showed that the composition $\Gcal\otimes \Gcal\to f^*\Acal''_{-2}$ corresponds to a morphism $f_+\Gcal\otimes f_+\Gcal \to \Acal.$ On the other hand, the composition corresponds  to the composition 
    $$\T(\Gcal\otimes \Gcal)\to \Acal_{-2}\to f^*\Acal''_{-2}$$
in $\Alg(\QCoh(Y)).$ This corresponds to the composition
$$f_\sharp\T(\Gcal\otimes \Gcal)\to f_\sharp\Acal_{-2}\to \Acal''_{-2}$$ in $\Alg_S.$ As $f_\sharp\T(\Gcal\otimes \Gcal)\simeq \T(f_+(\Gcal\otimes \Gcal))\simeq \T(f_+\Gcal\otimes f_+\Gcal)$, this corresponds to the composition
$$f_+\Gcal\otimes f_+\Gcal \xrightarrow{\psi'_{-2}}f_\sharp \Acal_{-2}\xrightarrow{u'_{-2}}\Acal''_{-2}.$$
Thus together with the strong assumption $(1)$ we have the corresponding commutative square
$$\xymatrix{
f_+\Gcal \otimes f_+\Gcal \ar[r]^-{\psi'_{-2}} \ar[d]_{u'_{-1}\otimes u'_{-1}}& f_\sharp\Acal_{-2} \ar[d]^{u'_{-2}}\\
\Gcal'' \otimes \Gcal'' \ar[r]_-{\psi''_{-2}}  &\Acal''_{-2}.
    }$$
\item We want to show that  $f^*\psi''_{-1}\circ u$ corresponds to $\psi''_{-1}\circ u'$ through the bijection $\Prm_{Y,\fin}(\Fcal,f^*(-))\xrightarrow{\simeq} \Prm_{S,\fin}(f_+\Fcal,-).$ Explicitly, this bijection sends a homogeneous polynomial law $\Gamma^i(\Fcal)\to f^*\Hcal$ to the compostion 
$$\Gamma^i(f_+\Fcal)\xleftarrow{\simeq}f_+\Gamma^i(\Fcal)\to f_+f^*\Hcal\xrightarrow{\epsilon_{\Hcal}}\Hcal.$$
We now prove that an inverse of $\Prm_{Y,\fin}(\Fcal,f^*(-))\xrightarrow{\simeq} \Prm_{S,\fin}(f_+\Fcal,-)$ could be given as follows: Any polynomial law $f_+\Fcal\leadsto \Hcal$ gives rise to a polynomial law $f^*f_+\Fcal\leadsto f^*\Hcal$ by pullback and precomposition with $\epsilon_{\Fcal}: \Fcal\to f^*f_+\Fcal$ gives a polynomial law $\Fcal\leadsto f^*\Hcal.$ Let $\psi:\Fcal\leadsto f^*\Hcal$ be a homgeneous component law of degree $i$ then the composition $\Prm_{Y,\fin}(\Fcal,f^*(-))\to \Prm_{S,\fin}(f_+\Fcal,-)\to \Prm_{Y,\fin}(\Fcal,f^*(-))$ sends $\psi$ to the composition
$$\Gamma^i(\Fcal)\xrightarrow{\Gamma^i(\eta_\Fcal)}\Gamma^i(f^*f_+\Fcal)\to f^*\Gamma^i(f_+\Fcal)\xleftarrow{\simeq}f^*f_+\Gamma^i(\Fcal)\to f^*f_+f^*\Hcal\xrightarrow{f^*\epsilon_{\Hcal}}f^*\Hcal.$$
Hence if suffices to show that the linear map $\Gamma^i(\Fcal)\to\Hcal$ induced by $\psi$ make the following diagram commute
$$
\xymatrix{
\Gamma^i(\Fcal)\ar[d]_{\Gamma^i(\eta_\Fcal)} \ar[r]^{\psi}& f^*\Hcal\\
\Gamma^i(f^*f_+\Fcal)\ar[d] & f^*f_+f^*\Hcal \ar[u]^{f^*\epsilon_\Hcal}\\
f^*\Gamma^i(f_+\Fcal)& f^*f_+\Gamma^i(\Fcal)\ar[l]^{\simeq} \ar[u]
}
$$
(note that the bottom map needs to be inverted so that the diagram make sense). We observe that it suffices to show that the following diagram 
$$
\xymatrix{
\Gamma^i(\Fcal)\ar[rdd]^{\eta_{\Gamma^i(\Fcal)}}\ar[d]_{\Gamma^i(\eta_\Fcal)} &\\
\Gamma^i(f^*f_+\Fcal)\ar[d] & \\
f^*\Gamma^i(f_+\Fcal)& f^*f_+\Gamma^i(\Fcal)\ar[l]^{\simeq} 
}
$$
commute as the diagram
$$
\xymatrix{
\Gamma^i(\Fcal)\ar[rdd]_{\eta_{\Gamma^i(\Fcal)}} \ar[r]^{\psi}& f^*\Hcal\\
 & f^*f_+f^*\Hcal \ar[u]^{f^*\epsilon_\Hcal}\\
& f^*f_+\Gamma^i(\Fcal) \ar[u]
}
$$
already commutes essentially by the definition of the unit $\eta_{\Gamma^i(\Fcal)}$, the counit $\epsilon_{\Hcal}$ and the unit $\eta_{\Hcal}$. Unwinding the definition of $f_+\Gamma^i(\Fcal)\to \Gamma^i(f_+\Fcal)$, we now need to show the following diagram
$$
\xymatrix{
\Gamma^i(\Fcal)\ar[ddrrr]^{\eta_{\Gamma^i(\Fcal)}}\ar[d]_{\Gamma^i(\eta_\Fcal)} &&& \\
\Gamma^i(f^*f_+\Fcal)\ar[d] & &&\\
f^*\Gamma^i(f_+\Fcal)& f^*f_+f^*\Gamma^i(f_+\Fcal) \ar[l] &f^*f_+\Gamma^i(f^*f_+\Fcal)\ar[l]&f^*f_+\Gamma^i(\Fcal)\ar[l] 
}
$$
commute. As we have the following small commutative square inside:
$$
\xymatrix{
\Gamma^i(\Fcal)\ar[ddrrr]^{\eta_{\Gamma^i(\Fcal)}}\ar[d]_{\Gamma^i(\eta_\Fcal)} &&& \\
\Gamma^i(f^*f_+\Fcal)\ar[drr]_{\eta_{\Gamma^i(f^*f_+\Fcal)}} & &&\\
&  &f^*f_+\Gamma^i(f^*f_+\Fcal)&f^*f_+\Gamma^i(\Fcal),\ar[l] 
}
$$
it suffices to show that the following diagram
$$
\xymatrix{
\Gamma^i(f^*f_+\Fcal)\ar[d]\ar[drr]^{\eta_{\Gamma^i(f^*f_+\Fcal)}} & &&\\
f^*\Gamma^i(f_+\Fcal)& f^*f_+f^*\Gamma^i(f_+\Fcal) \ar[l] &f^*f_+\Gamma^i(f^*f_+\Fcal)\ar[l]}
$$
commute, which is true essentially by the definition of the unit $\eta_{\Gamma^i(f^*f_+\Fcal)}$, the counit $\epsilon_{\Gamma^i(f_+\Fcal)}$ and the unit $\eta_{\Gamma^i(f_+\Fcal)}$.\\
Return to our problem showing that $f^*\psi''_{-1}\circ u$ corresponds to $\psi''_{-1}\circ u'$ through the bijection $\Prm_{Y,\fin}(\Fcal,f^*(-))\xrightarrow{\simeq} \Prm_{S,\fin}(f_+\Fcal,-)$, it suffices to show that $$(f^*(\psi''_{-1}\circ u'))\circ \eta_{\Fcal}=f^*\psi''_{-1} \circ u,$$ which is true as the left hand side is 
$$f^*(\psi''_{-1})\circ( f^*u'\circ \eta_{\Fcal})=f^*\psi''_{-1} \circ u.$$
\item Finally, we show that $u_i\circ \psi_i$ corresponds to $u'_i\circ \psi'_i$ through the bijection $\Prm_{S,\fin}(f_+\Fcal,\Acal''_i)\simeq \Prm_{Y,\fin}(\Fcal,f^*\Acal''_i).$ Indeed, we have $$f^*(u'_i\circ \psi'_i)\circ \eta_\Fcal=f^*u'_i\circ (f^*\psi'_i\circ \eta_\Fcal)=f^*u'_i\circ (\eta_\Acal\circ \psi_i)=(f^*u'_i\circ \eta_\Fcal)\circ \psi_i=u_i\circ \psi_i.$$
\end{itemize}
Thus the tuples $(u,u_{-2},\dots,u_d)$ correspond to the tuples $(u',u'_{-2},\dots,u_d)$ making the following diagrams
$$\xymatrix{
f_+\Gcal \otimes f_+\Gcal \ar[r]^-{\psi'_{-2}} \ar[d]_{u'_{-1}\otimes u'_{-1}}& f_\sharp\Acal_{-2} \ar[d]^{u'_{-2}}\\
\Gcal'' \otimes \Gcal'' \ar[r]_-{\psi''_{-2}}  &\Acal''_{-2}
    },
    \xymatrix{
     f_+\Gcal\ar[d]_{u'_{-1}} &f_+\Fcal\ar[l]_{\psi'_{-1}} \ar[d]^{u'}\ar[r]^{\psi'_{i}}& f_\sharp\Acal_i\ar[d]^{u'_{i}}\\ \Gcal''&\Fcal''\ar[l]_{\psi''_{-1}}\ar[r]^-{\psi''_{i}} &\Acal''_i
    }$$
commute and we are done.
\end{proof} 
\end{theorem}
Our main application of the theorem is the following:
\begin{lemma}\label{ex:condition}
Let $\Fcal$ be a finite locally free $\Ocal_S$-module and let $$p:\Pbb(\Fcal)=\Proj(\Sym(\Fcal^\vee))\to S$$ be the relative proj construction. Then the line bundles $\Ocal_{\Pbb(\Fcal)}(-1)$ satisfy the assumption canonical morphism $p_+(\Gamma^i(\Ocal(-1)))\to \Gamma^i(p_+\Ocal(-1))$ of Theorem \ref{thm:main_compute}. 
\begin{proof}
By Remark \ref{rmk:module_lef}, for $j\geq 0,$ $$p_+(\Ocal(-j))\simeq \Hrm_0(\R p_*(\Ocal(-j)^\vee))^\vee.$$
As $\Ocal(-j)$ is flat, we have $\Ocal(-j)^\vee=\R\Hom_{\Ocal_{\Pbb(\Fcal)}}(\Ocal(-j),\Ocal_{\Pbb(\Fcal)})\simeq \Ocal(j)$ and hence $p_{+}\Ocal(-j)\simeq p_*(\Ocal(j))^\vee \simeq \Sym^j(\Fcal^\vee)^\vee.$ Note that these choices of isomorphisms comes from the identification $p_*\Ocal(j)\simeq \Sym^j(\Fcal^\vee)$ and they are canonical as part of relative proj construction (see \cite[\href{https://stacks.math.columbia.edu/tag/01OC}{Tag 01OC}]{stack}). We are done if the following diagram
$$
\xymatrix{
p_+(\Gamma^i(\Ocal(-1)))\ar[d] \ar[r]& \Gamma^i(p_+\Ocal(-1)) \ar[dl]\\
\Sym^i(\Fcal^\vee)^\vee
}
$$
commutes. Indeed, by adjunction we consider the following diagram
$$
\xymatrix{
\Gamma^i(\Ocal(-1))\ar[d] \ar[r]& p^*\Gamma^i(p_+\Ocal(-1)) \ar[dl]\\
p^*(\Sym^i(\Fcal^\vee)^\vee).
}
$$
Now it suffices to show that the choice of isomorphism $p_*\Ocal(1)\simeq \Fcal^\vee$ gives rise to the morphism $\Gamma^i(\Ocal(-1))\to p^*(\Sym^i(\Fcal^\vee)^\vee).$ This is a well-know property of relative proj construction that the the canonical morphism $\Sym(\Fcal^\vee)\to p_*(\oplus_{n\geq 0} \Ocal(n))$ is graded (see \cite[\href{https://stacks.math.columbia.edu/tag/01NR}{Tag 01NR}, \href{https://stacks.math.columbia.edu/tag/01OC}{Tag 01OC}]{stack}).
\end{proof}
\end{lemma}
\begin{construction}\label{con:restrict}
Consider the setting of Lemma \ref{ex:condition}. Let $\Psi\in \poly^d(S)$ given by:
\begin{itemize}
    \item $\psi_{-1}:\Fcal \to \Gamma^m(\Fcal)$ is the canonical polynomial law $\lrm_\Fcal^m.$
    \item $\psi_i:\Fcal\to \Ocal_S$ be a homogeneous polynomial law of degree $im$ for $i=0,1,\dots,d.$
\end{itemize}
Let $\psi|_{\Ocal(-1)}$ denote the composition $p^*\psi \circ (\Ocal(-1)\to p^*\Fcal).$ This restriction together with the law $\lrm_{\Ocal(-1)}^m$ give a weighted polynomial law $\Psi|_{\Ocal(-1)}\in \poly^d(\Pbb(\Fcal))$ given by $$\Psi|_{\Ocal(-1)}=(\Ocal(-1),\Gamma^m(\Ocal(-1)),\Ocal_{\Pbb(\Fcal)},\dots,\Ocal_{\Pbb(\Fcal)},\lrm^m_{\Ocal(-1)},(\psi_i|_{\Ocal(-1)})_{0\leq i\leq d}).$$ 
\end{construction}
\begin{corollary}\label{cor:main2}
Consider the setting of Construction \ref{con:restrict}, we have $\Psi|_{\Ocal(-1)}$ is a lift of $\Psi$ along $p.$ As a consequence of Proposition \ref{prop:lift_alg}, we obtain $p_\sharp \Cl(\Psi|_{\Ocal(-1)})\simeq \Cl(\Psi).$ 
\begin{proof}
By Lemma \ref{ex:condition} and Theorem \ref{thm:main_compute}, $\Psi|_{\Ocal(-1)}$ is a lift along $p$ of a weighted polynomial law $\Psi_{S}\in \poly^d(S)$. Moreover, the proof of Theorem \ref{thm:main_compute} gives the formula of $\Psi_S$ which consists of the following weights:
\begin{itemize}
    \item The polynomial law $p_+\Ocal(-1)\leadsto p_+\Gamma^m(\Ocal(-1))$ given by applying $p_+$ to the identity $\Gamma^m(\Ocal(-1))\to \Gamma^m(\Ocal(-1))$ and takes the composition with the isomorphsim $$\Gamma^m(p_+\Ocal(-1))\to p_+\Gamma^m(\Ocal(-1)).$$ In the proof of Lemma \ref{ex:condition}, we know that $p_+$ acts as $p_*((-)^\vee)^\vee$ on line bundles $\Ocal(j)$ so we obtain $\id_{\Gamma^m(\Fcal)}.$
    \item The polynomial law $p_+(\Ocal(-1))\leadsto \Ocal_S$ given by applying $p_+$ to the morphism $\Ocal(-i)\to \Ocal_{\Pbb(\Fcal)}$ and taking the composition with the isomorphism $$\Gamma^i(p_+\Ocal(-1))\to p_+\Gamma^i(\Ocal(-1)).$$ As we know that $p_+$ acts as $p_*((-)^\vee)^\vee$ on line bundles $\Ocal(j)$, the result polynomial law is $\psi_i$.\hfil\qedhere
\end{itemize}
\end{proof}
\end{corollary}
\begin{remark}\label{rmk:com}
In Proposition \ref{prop:lift_alg}, if it occurs that $\Cl(\Psi_Y)$ is commutative then we can identify it with a scheme over $S$ by taking its relative spectrum over $Y.$ For example, if a weighted polynomial law $\Psi$ satisfies $\Gcal$ is a line bundle and all sheaves $\Acal_i$ are commutative then the explicit construction implies that $\Cl(\Psi)$ is commutative.
In particular, this is true for polynomial law $\Psi|_{\Ocal(-1)}$ in Construction \ref{con:restrict}. As we will see in the next section, these schemes are exactly the varieties appearing in \cite{Bergh_1987} or more recently by \cite{Chapman_2015} which are all generalized in \cite[Appendix A]{KL_2015}.
\end{remark}
\begin{remark}
We can construct rank $2$ vector bundles that satisfy the assumption of Theorem \ref{thm:main_compute} by observing that if $\Lcal$ is such a line bundle over $Y$ then so is $\Lcal\oplus \Ocal_Y.$ For we have
$$\Gamma^i(\Lcal\oplus \Ocal_Y)\simeq \Sym^i((\Lcal\oplus \Ocal_Y)^\vee)^\vee \simeq (\oplus_{k=0}^i\Lcal^{\otimes -i})^\vee \simeq \oplus_{k=0}^i \Lcal^{\otimes i}\simeq \oplus_{k=0}^i \Gamma^i(\Lcal)$$
and hence the natural map $f_+\Gamma^i(\Lcal\oplus \Ocal_Y)\to \Gamma^i(f_+(\Lcal\oplus \Ocal_Y))$ is an isomorphism.
\end{remark}

\subsection{Comparison with Krashen and Lieblich's explicit construction}
In this section, we show that how to deduce Krashen and Lieblich's explicit construction of Roby's Clifford algebra given in \cite[Proposition A.1.1]{KL_2015} from the above theory. We start by recalling some content in \cite{KL_2015}. For a morphism of schemes $\phi:X\to Y$ over a scheme $S$, \textit{the Clifford functor of $\phi$} is the functor
$$\Ccal\Fcal_\phi:\mathrm{LFAlg}/\Ocal_S\to \Set$$
given by the formula
$$\Ccal\Fcal_\phi(\Acal):=\Hom_{\Alg_Y}(\phi_*\Ocal_X,\pi_Y^*\Acal),$$
where $\mathrm{LFAlg}/\Ocal_S$ is the category of sheaves of algebras over $S$ that are locally free sheaves of $\Ocal_S$-modules, and $\Alg_Y$ is the category of quasi-coherent sheaves of $\Ocal_Y$-algebras. The quasi-coherent sheaf of algebras over $S$ that corepresents the Clifford functor is called \textit{the Clifford algebra of $\phi$}.

To show that the above construction is a generalization of various definitions of Clifford algebras, Krashen and Lieblich pick a morphism $\phi$ corresponding to each definition of Clifford algebra so that that Clifford algebra is exactly the Clifford algebra of $\phi$. This is done in  \cite[Appendix A]{KL_2015} which we recalled in Example \ref{ex:KL}. They consider Roby's Clifford algebra in Section A.1 of Appendix A, as well as Clifford algebra definitions of other authors. Then in Section A.4, they give a new common generalization of these algebras by showing that this new one is the Clifford algebra of some morphism $\phi$.
\begin{theorem}[{\cite[Proposition A.4.1]{KL_2015}}]\label{thm:KL}
Let $X$ be the hypersurface in the weighted projective space $\Pbb_{m,1,\dots,1}$ defined by the degree $md$ homogeneous equation $$x_0^d=x_0^{d-1}+x_0^{d-2}f_{2m}+\cdots+f_{dm}.$$
Let $X\to \Pbb^{n-1}_{k}$ be the degree $d$ morphism given by dropping the $x_0$-coordinate. Then $\Ccal\Fcal_{\phi}$ is corepresented by the algebra $\Cl(f_m,\cdots,f_{dm}).$ 
\end{theorem}
We now describe a different proof of this theorem, slightly generalizing it on the way: The base $k$ is allowed to be an arbitrary scheme, and the projective space can be that of an arbitrary vector bundle.  The precise statement is in Corollary \ref{cor:final}

Let $\Psi\in \poly^d(S)$ be as in Construction \ref{con:restrict} with the additional assumption that $\psi_0=1$. By  Remark \ref{rmk:com}, the sheaf $\Cl(\Psi|_{\Ocal(-1)})$ is commutative and has a relative spectrum over $\Pbb(\Fcal).$ We give a concrete description of this relative spectrum in terms of functor of points. Recall that $\Proj_S(\Acal)$ of a quasi-coherent graded $\Ocal$-algebra $\Acal$ sends $(A,x:\Spec(A)\to S)$ to the set of equivalent classes
$$
\left\{[\text{surjective }\Acal(A,x)\to L_A]\bigg| \begin{matrix}
L_A \text{ is a line bundle over }A \text{ s.t } \\ \Acal(A,x)^{(t)}\to \oplus_{j\geq 0}L^{\otimes j} \text{ is graded for some }t\geq 1    
\end{matrix}
\right\}.
$$
The functor of points of $\Pbb(\Fcal)$ is then the $S$-functor that sends $(A,x:\Spec(A)\to S)$ to the set of equivalent classes $$
\left\{[\text{surjective }\Fcal(A,x)^\vee\to L_A]
\right\}.
$$
We let $\Pbb_m(\Fcal\oplus \Ocal_S)$ denote $\Proj(\Sym(\Fcal^\vee\oplus \Ocal_S))$ and using the formula $\Sym(\Fcal^\vee\oplus \Ocal_S)\simeq \Sym(\Fcal^\vee)\otimes \Sym(\Ocal_S),$ the grading of $\Sym(\Fcal^\vee\oplus \Ocal_S)$ is obtained by imposing the grade $m$ on $\Ocal_S$ (and grade $1$ on $\Fcal$). From now we freely use these gradings. The functor points of $\Pbb_{m}(\Fcal\oplus \Ocal_S)$ is the $S$-functor that sends $(A,x)$ to
$$
\left\{[\text{surjective }\Sym^m(\Fcal(A,x)^\vee)\oplus A^\vee\to L_A]
\right\}.
$$
For each pair $(A,x)$, we can identify polynomial laws $\psi_i|_{(A,x)}:\Gamma^{mi}(\Fcal(A,x))\to A$ as elements of $\Sym^{mi}(\Fcal(A,x)^\vee)$ by the isomorphism $\Gamma^{mi}(\Fcal(A,x))\simeq \Sym^{mi}(\Fcal(A,x)^\vee)^\vee.$ These elements give an element 
$$\psi_{(A,x)}=\psi_d|_{(A,x)}\oplus\psi_{d-1}|_{(A,x)}\oplus\cdots\oplus\psi_1|_{(A,x)}+\id_A$$
in $$\Sym^{md}(\Fcal(A,x)^\vee)\oplus \Sym^{m(d-1)}(\Fcal(A,x)^\vee)\oplus\cdots \oplus \Sym^m(\Fcal(A,x)^\vee)\oplus A^\vee.$$
Define $X=X_\Psi$ the subfunctor of $\Pbb(\Fcal\oplus \Ocal_S)$ that sends $(A,x)$ to
\begin{equation}\label{eq:X_psi}
\left\{[\text{surjective }\Sym^m(\Fcal(A,x)^\vee)\oplus A^\vee\to L_A]\mid \psi_{(A,x)}\mapsto 0 \in L_A^{\otimes d}
\right\}.
\end{equation}
We now define a projection morphism $X\to \Pbb(\Fcal):$
\begin{lemma}\label{lm:well_def}
If $[\Sym^m(\Fcal(A,x)^\vee)\oplus A^\vee\to L_A]\in X(A,x)$ then the composition $$\Sym^m(\Fcal(A,x)^\vee)\to\Sym^m(\Fcal(A,x)^\vee)\oplus A^\vee\to L_A$$
is surjective. Consequently, we obtain a natural transformation $\phi:X\to \Pbb(\Fcal).$
\end{lemma}
\begin{proof}
By taking localization, we can assume that $A$ is a local ring, and $L_A$ and $\Fcal(A,x)$ are free. Let $\{e_1,\dots,e_n\}$ be a basis of $\Fcal(A,x),$ then we can write $\psi_i$ as a homogeneous polynomial $\psi_i(e_1^*,\dots,e_n^*)$ of degree $mi$ in the dual basis $\{e_1^*,\dots,e_n^*\}$. We also let $\{e_0\}$ be the basis of $A$ corresponding to $1_A$, and consider $\id_A:A\to A$ as the basis $\{e_0^*\}.$ Let $h$ denote the corresponding map 
$\Sym^m(\Fcal(A,x)^\vee)\oplus A^\vee \to L_A^\vee.$ The condition that $\psi_{(A,x)}$ is mapped to $0$ becomes
\begin{equation}\label{eq:psi-d}
-\psi_d(h_d(e^*))+\psi_{d-1}(h_{d-1}(e^*))h(e_0^*)+\cdots+\psi_m(h_{m}(e^*))(h(e_0^*))^{d-1}+(h(e_0^*))^{d}=0
\end{equation}
in $L_A^{\otimes d},$ where we let $h_j(e^*)$ denote any tuple of $j$ elements in $$\{h((e_{1}^*)^{i_1}\cdots(e_n^*)^{i_n})\mid i_1+\cdots +i_n = m\}.$$ Let $\{e\}$ be a basis of $L_A.$ Since $h$ is surjective, there exists $a_0\in A$ and $a_{(i_1,i_2,\dots,i_n)}\in A$ such that $$e^*=a_0h(e_0^*)+\sum_{i_1+\cdots+i_n=m} a_{(i_1,\dots,i_n)}h((e_{1}^*)^{i_1}\cdots(e_n^*)^{i_n}).$$ 
Writing 
$h((e_{1}^*)^{i_1}\cdots(e_n^*)^{i_n})=b_{(i_1,\dots,i_n)}e^*$, this is equivalent to 
\begin{equation}\label{eq:unit}
1=a_0b_0+\sum_{i_1+\cdots+i_n=m}a_{(i_1,i_2,\dots,i_n)}b_{(i_1,\dots,i_n)}
\end{equation}
in $A$. On the other hand, the equation \eqref{eq:psi-d} reduces to
$$-\psi_d(b_{(i_1,\dots,i_n)}e^*)+\psi_{d-1}(b_{(i_1,\dots,i_n)}e^*)b_0e^*+\cdots+\psi_m(b_{(i_1,\dots,i_n)}e^*)b_0^{d-1}(e^*)^{d-1}+b_0^d(e^*)^d=0.$$
In particular, $b_0^d$ can be written as a homogeneous polynomial in $b_0,b_{(i_1,\dots,i_n)}$ where $b_0$ is of smaller degrees. Hence any prime ideal that contains $b_{(i_1,\dots,i_n)}$ also contains $b_0^d$, and therefore $b_0$. By \eqref{eq:unit}, such prime ideal contains $1$ so in fact, $$(b_{i_1,\dots,i_n}\mid i_1+\cdots+i_n=m)=A.$$ 
It follows  that $h|_{\Sym^m(\Fcal(A,x))}$ is surjective. More explicitly, there are $a'_{(i_1,\dots,i_n)}\in A$ such that $1=\sum_{i_1+\cdots+i_m}a'_{(i_1,\dots,i_n)}b_{(i_1,\dots,i_n)}$ and $h$ maps $\sum_{i_1+\cdots+i_m}a'_{(i_1,\dots,i_n)}(e_{1}^*)^{i_1}\cdots(e_n^*)^{i_n}$ to $e^*.$
\end{proof}
\begin{theorem}\label{thm:main3}
The relative spectrum $\Spec(\Cl(\Psi|_{\Ocal(-1)})\to \Pbb(\Fcal)$ is isomorphic to $\phi:X_\Psi\to \Pbb(\Fcal)$.
\end{theorem}
\begin{proof}
We are going to show that $\phi_*\Ocal_{X_\Psi}\simeq \Cl(\Psi|_{\Ocal(-1)})$
as sheaves of algebras over $\Pbb(\Fcal).$ For each point $(A,x)\in \CAlg(\Pbb(\Fcal)),$ 
$$\phi_*\Ocal_{X_\Psi}(A,x)=\Ocal_{X_\Psi}(A',x')=A',$$
where $A'$ fits into the following pullback diagram
$$
\xymatrix{
\Spec(A') \ar[r] \ar[d]& \Spec(A)\ar[d]^{x}\\
X \ar[r] & \Pbb(\Fcal).
}
$$
On the other hand, by Lemma \ref{lm:func_cliff}, the sheaf $\Cl(\Psi|_{\Ocal(-1)})$ corepresents the functor which assigns to each $(A,x)\in \CAlg(\Pbb(\Fcal))$ the algebra
$\Cl((\Psi|_{\Ocal(-1)})|_{(A,x)}).$ Thus it suffices to show that there exists an isomorphism $\Spec(A')\simeq \Spec(\Cl((\Psi|_{\Ocal(-1)})|_{(A,x)}))$ as a natural isomorphism of functors over $\CAlg(S).$ Let $[\Fcal(A,x')^\vee\to L_A]$ be a representative of $(A,x)$

We first give a description for the functor of points of $\Spec(A').$ For each $(B,y)\in \CAlg(S)$, we have the following pullback diagram in $\Set$:
$$
\xymatrix{
\Spec(A')(B,y) \ar[r] \ar[d]& \Spec(A)(B,y)\ar[d]\\
X(B,y) \ar[r] & \Pbb(\Fcal)(B,y),
}
$$
and hence one can compute $\Spec(A')(B,y)$ as $\Spec(A)(B,y)\times_{\Pbb(\Fcal)(B,y)} X_\Psi(B,y)$. We have 
$$\Spec(A)(B,y)=\{f:A\to B\mid S(f)(x')=y\}$$
where $x'=(\Pbb(\Fcal)\to S)\circ x$. Let $[\Fcal(B,y)^\vee\to L_B]$ be a representative of $(B,y)$. Unwinding the equation $S(f)(x')=y$, we obtain isomorphisms making the following diagram 
$$
\xymatrix{
\Fcal(A,x')^\vee\otimes_A B \ar[r] \ar[d]_{\simeq}& L_A\otimes_A B \ar[d]^{\simeq}\\
\Fcal(B,y)^\vee \ar[r] & L_B
}
$$
commute. On the other hand, each $[\Sym^m(\Fcal(B,y)^\vee)\oplus B^\vee\to L'_B]\in X_\Psi(B,y)$ maps to a class $[\Fcal(B,y)^\vee\to L_B]$ in $\Pbb(\Fcal)(B,y)$ satisfies there exists an isomorphism $L_B \to L'^{\otimes m}_B$ making the following diagram commute
$$
\xymatrix{
\Sym^m(\Fcal(B,y)^\vee) \ar[d]_{=}\ar[r] & L'_B \ar[d]\\
\Sym^m(\Fcal(B,y)^\vee)  \ar[r]_-{\Sym^m} & L^{\otimes m}_B.
}
$$
Thus, the set $\Spec(A)(B,y)\times_{\Pbb(\Fcal)(B,y)} X_\Psi(B,y)$ is the set of pairs $$(f, B^\vee\to (L_A\otimes B)^{\otimes m})$$
satisfying $\psi_{(B,y)}$ is mapped to $0$ in $L'_B$. 

We now give the description of functor of points of $\Spec(\Cl((\Psi|_{\Ocal(-1)}))).$ Recall that each weight $\psi_i|_{\Ocal(-1)}$ is the composition law 
$$\Ocal(-1)\to p^*\Fcal\overset{p^*\psi_i}{\leadsto} \Ocal_{\Pbb(\Fcal)},$$
so its restriction $(\psi_i|_{\Ocal(-1)})|_{(A,x)}$ is the composition of polynomial laws
$$L_A^\vee \xhookrightarrow{} \Fcal(A,x')\overset{\psi_i|_{(A,x')}}{\leadsto} A$$
which we denote by $\psi_i|_{L_A^\vee}.$ It follows that $\Spec\left(\Cl((\Psi|_{\Ocal(-1)})|_{(A,x)})\right)(B,y)$ is the set of ring homomorphisms $$g:\Cl((\Psi|_{\Ocal(-1)})|_{(A,x)})\to B$$
such that 
$$
S(g)(x'\circ (\Spec(\Cl((\Psi|_{\Ocal(-1)})|_{(A,x)}))\to \Spec(A))) = y.
$$
Such a homomorphism $g$ is equivalent to a pair $((L_A^\vee)^{\otimes m}\to B,f:A\to B)$ satisfying the polynomial law
$$\phi\circ (\psi_d|_{L_A^\vee})+(L_A^\vee \to B)(\phi\circ \psi_{d-1}|_{L_A^\vee})+\cdots+\id_B^{d}\circ (L_A^\vee \to B)=0,$$
where $L_A^\vee\to B$ is given by the composition $L_A^\vee \xrightarrow{z\mapsto z^{\otimes m}} (L_A^\vee)^{\otimes m}\to B.$ 

Comparing these two descriptions of $\Spec(A')(B,y)$ and $\Spec(\Cl(\Psi|_{\Ocal(-1)})|_{(A,x)})(B,y)$, we concludes that 
$$\Spec(A') \simeq \Spec(\Cl((\Psi|_{\Ocal(-1)})|_{(A,x)})).            \hfil\qedhere$$
\end{proof}
\begin{corollary}\label{cor:final}
Let $\phi:X_\Psi\to \Pbb(\Fcal)$ the projection map. Then the Clifford functor $\Ccal\Fcal_\phi:\QCoh(S)\to \Set$ is corepresented by the Clifford algebra $\Cl(\Psi)$.
\begin{proof}
This is true since
\begin{align*}
\Ccal\Fcal_\phi(\Acal) &= \Hom_{\Alg_{\Pbb(\Fcal)}}(\phi_*\Ocal_{X_\Psi},p^* \Acal)\\
&\simeq \Hom_{\Alg(\Pbb(\Fcal))}(\Cl(\Psi|_{\Ocal(-1)}),p^*\Acal) \quad (\text{by Theorem }\ref{thm:main3})\\
&\simeq \Hom_{\Alg_S}(p_\sharp\Cl(\Psi|_{\Ocal(-1)}),\Acal) \quad (\text{by Theorem }\ref{ad_module})\\
&\simeq \Hom_{\Alg_S}(\Cl(\Psi),\Acal) \quad (\text{by Corollary }\ref{cor:main2}).\hfil\qedhere
\end{align*}
\end{proof}
\end{corollary}
\section{Derived Clifford algebras}
As an application of the  realization of Clifford algebra construction as a left adjoint, we define a derived version of Clifford algebra with noncommutative coefficients. We start with the case of a commutative ring $k,$ which has some unique features in its relation with classical Clifford algebra. Later we will extend our constructions to any animated commutative ring. 
\subsection{Derived Clifford algebras over commutative rings}
We let $\Mod_k^\an$ denote the animation of $\Mod_k$, which is also the category of connective objects of derived category $\Dcal(R)$ (\cite[5.1.7]{C-S_2019}).

By applying \cite[Corollary 4.7.3.18]{HA} to the adjunction $(\Fr\dashv i):\Alg_k\leftrightarrows \Set$, the category $\Alg_k$ is generated under colimits by $\Alg_k^\sfp$ where $\Alg_k^\sfp$ is the full subcategory of retracts of free algebras $k\langle X_1,\dots, X_n\rangle.$ We let $\Alg_k^\an$ denote the animation of $\Alg_k$, which by \cite[5.1.4]{C-S_2019}, is the sifted cocompletion of $\Alg_k^{\sfp}$ (see \cite[\href{https://kerodon.net/tag/04BJ}{Tag 04BJ}]{ker} for the definition of cocompletion). By \cite[Proposition 7.1.4.18]{HA}, the $\infty$-category $\Alg^\an_k$ is equivalent to $\Alg_k^{(1),\mathrm{cn}}$ of connective $\Ebb_1$-algebras over $k$ as well as the $\infty$-localization of the category of simplicial associative algebras over $R$ at weak equivalences. 

In  this section, we consider a special variant of Definition \ref{def:w_poly}. For any commutative ring $R$, we define the category $\poly_{h,d}(k)$ as the comma category $(\Gamma^d_k \downarrow \id_{\Mod_k}).$ The Clifford algebra functor $\poly_{h,d}(k)\to \Alg_k$ is the left adjoint of functor $$F_{h,d}:\Alg_k\to \poly_{h,d}(k), A\mapsto (\Gamma^d(A)\to A, z^{[d]}\mapsto z^d).$$ A natural derived version of $\poly_{h,d}(k)$ is then constructed as follows. 
\begin{definition}
We define the $\infty$-category of \textit{derived homogeneous polynomial forms of degree $d$ with noncommutative coefficients}, denoted by $\drpoly_{h,d}(k)$, as the oriented fiber product $\Mod_k^\an \times^{\to}_{\Mod_k^\an} \Alg_k^\an$ where we use  $\Lrm\Gamma^d:\Mod_k^\an\to\Mod_k^\an$ and the forgetful functor $\Alg_k^\an\to \Mod_k^\an$ in the construction (the definition of oriented fiber product is in \cite[\href{https://kerodon.net/tag/01KF}{Tag 01KF}]{ker}).
\end{definition}
\begin{remark}
The objects of $\drpoly_{h,d}(k)$ are the triples $(K,A,\psi)$ where $K \in \Mod_k^\an,A\in \Alg_k^\an,$ and $\psi$ is a morphism from $\Lrm\Gamma^d(K)\to A$ in $\Mod_k^\an$.
\end{remark}
On the other hand, one can also think about animation of $\poly_{h,d}(k)$ as a candidate of its derived version. We are going to show that these two constructions are equivalent. To obtain animation of $\poly_{h,d}(k)$, we need to show that $\poly_{h,d}(k)$ is cocomplete and generated under colimits by $\poly_{h,d}(k)^\sfp.$ 
\begin{lemma}\label{lm:plus_adjunction}
Let $F:\Ecal \to \Dcal, G: \Ccal\to \Dcal$ be functors between ordinary categories and assume that $\Ccal$ has finite coproducts. If $G$ admits a left adjoint then the projection functor $(F\downarrow G) \to \Ecal\times \Ccal, (x,y,Fx\to Gy)\mapsto (x,y)$ admits a left adjoint. 
\end{lemma}
\begin{proof}
Let $H$ be the left adjoint of $G$. The functor $$\Ecal \times \Ccal \to (F\downarrow G),\ (x,y)\mapsto (x,HFx\sqcup y,Fx\xrightarrow{\eta_{Fx}} GHFx\rightarrow G(HFx  \sqcup y))$$ will does the job where $\eta$ is the unit transformation of the adjunction $F\dashv G$ and $GHFx\to G(HFx\sqcup y)$ is induced by the inclusion $HFx\hookrightarrow HFx\sqcup y$. Indeed, for any $(x',y',f:Fx'\to Gy')$ in $(F\downarrow G),$ we have
$$\Hom_{ (F\downarrow G)}((x,HFx,Fx\to G(HFx\sqcup y)), (x',y',f:Fx'\to Gy'))$$ is bijective to pairs $(u,v)$ in $\Hom_{\Ecal}(x,x')\times \Hom_{\Ccal}(HFx\sqcup y,y')$ such that the following square 
$$
\xymatrix{
Fx \ar[r] \ar[d]_-{Fu} & G(HFx\sqcup y) \ar[d]^{Gv}\\
Fx' \ar[r]_-{f} & Gy'
}
$$
commutes. This square is equivalent to commutative square
$$
\xymatrix{
HFx \ar[r] \ar[d]_-{HFu} & HFx \sqcup y \ar[d]^{v}\\
HFx' \ar[r]_-{\epsilon_{y'}\circ Hf} & y',
}
$$
where the horizontal top arrow is the inclusion. Hence it is equivalent to a morphism $v':y\to y'$ in $\Ccal$ and the commutative square 
$$
\xymatrix{
HFx\ar[r]^-{\id_{Fx}}\ar[d]_{HFu} & HFx \ar[d]^{v|_{HFx}}\\
HFx'\ar[r]_-{\epsilon_{y'}\circ Hf}& y',
}
$$
which implies that $v$ is completely determined by $u$ and $v'$.
\end{proof}
Applying Lemma \ref{lm:plus_adjunction} to $\poly_{h,d}(k)=(\Gamma^d_k\downarrow i_{\Mod_k})$ as $i_{\Mod_k}$ has left adjoint the free algebra functor, we obtain a left adjoint $\Mod_k\times \Alg_k\to \poly_{h,d}(k)$. We now invoke Corollary 2.3 from \cite{ZM}, as stated below.
\begin{lemma}\label{lm:ani_adjoint}
Let $(F\dashv G): \Dcal \leftrightarrows \Ccal$ be an adjunction such that:
\begin{enumerate}
    \item $\Dcal$ admits filtered colimits and reflexive coequalizers and $G$ preserves filtered colimits and reflexive coequalizers.
    \item $\Ccal$ is cocomplete and generated under colimits by $\Ccal^\sfp$.
    \item $G$ is conservative.
\end{enumerate}
Then $\Dcal$ is cocomplete and generatd under colimits by $\Dcal^\sfp$ and we obtain an adjunction pair $(\Ani(F)\dashv \Ani(G)):\Ani(\Dcal)\leftrightarrows \Ani(\Ccal)$ of $\infty$-categories.
\end{lemma}
So we obtain a left adjoint of (animated) projection functor $\Ani(\poly_{h,d}(k))\to \Mod_k^\an\times \Alg_k^\an.$ We can choose $$\{(K,\T(\Gamma^d(K)))*A,K\hookrightarrow \T(\Gamma^d(K)))*A)\mid K \in \Mod_k^\sfp, A\in \Alg_k^\sfp\}$$ as a collection of compact projective generators of $\Ani(\poly_{h,d}(k)).$ We now prove the following higher version of Lemma \ref{lm:plus_adjunction}. Although the following lemma implies Lemma \ref{lm:plus_adjunction}, but its proof is not so explicit so we keep the proof of Lemma \ref{lm:plus_adjunction} for reader's convenience. We will use the following Lurie's notion of relative adjunctions from \cite[Section 7.3.2]{HA}. Consider the following commutative diagram of $\infty$-categories
$$
\xymatrix{
\Ccal\ar[dr]\ar[r]^{G} & \Dcal \ar[d]\\
&\Ecal
}
$$
where $\Dcal \to \Ecal$ and $\Ccal\to \Ecal$ are isofibrations (see \cite[\href{https://kerodon.net/tag/01ES}{Tag 01ES}]{stack} for the definition of isofibration, which is called categorical fibration in \cite{HTT}). We say that $G$ \textit{admits a left adjoint relative to $\Ecal$} if it admits a left adjoint $H$ and the unit map $x\to GHx$ is mapped to an equivalence in $\Ecal.$ It suffices for us to use the following lemma which is a direct consequence of \cite[Proposition 7.3.2.6]{HA}
\begin{lemma}\label{lm:relative_adj}
Suppose we have a commutative diagram of $\infty$-categories
$$
\xymatrix{
\Ccal\ar[dr]_-{r}\ar[r]^{G} & \Dcal \ar[d]^-{q}\\
&\Ecal
}
$$
where $\Ccal\to\Ecal,\Dcal\to \Ecal$ are Cartesian fibrations. If for every $z\in \Ecal,$ the induced functor $$G_z:\Ccal_z\to \Dcal_z$$ admits a left adjoint and $G$ sends $r$-cartesian morphisms to $q$-cartesian morphisms then $G$ admits a left adjoint relative to $\Ecal.$
\end{lemma}
\begin{proof}
We make a modification that our functors $r,q$ are cartesian fibrations instead of locally cartesian as in \cite[Proposition 7.3.2.6]{HA}. This gives us a sufficient condition to have a relative adjunctions (but not a necessary as in \cite[Proposition 7.3.2.6]{HA}). 
\end{proof}
\begin{lemma}\label{lm:plusdrived_adjunction}
Let $F:\Ecal \to \Dcal, G: \Ccal\to \Dcal$ be functors between $\infty$-categories and assume that $\Ccal$ has finite coproducts. If $G$ admits a left adjoint then the projection functor $(\Dcal \times^{\to}_{\Dcal} \Ccal) \to \Dcal\times \Ccal, (x,y,Fx\to Gy)\mapsto (x,y)$ admits a left adjoint.
\end{lemma}
\begin{proof}
Consider the following commutative diagram of $\infty$-categories
$$
\xymatrix{
\Fun(\Delta^1,\Dcal)\times_{\Fun(\{1\},\Dcal)}\Ccal \ar[r]^-{p}\ar[dr]_-{r} & \Dcal\times \Ccal\ar[d]^{q}\\
& \Dcal,
}
$$
where $p$ sends $(x\to Gy,y)$ to $(x,y)$, $q$ is the projection to the first factor, and $r$ sends $(x\to Gy,y)$ to $x.$ We want to show that $p$ is a relative adjunction by Lemma \ref{lm:relative_adj}.

By \cite[Corollary 2.4.7.12]{HTT} and \cite[\href{https://kerodon.net/tag/0478}{Tag 0478}]{ker}, the functors $r, q$ are cartesian fibrations. Moreover, the morphisms in $\Fun(\Delta^1,\Dcal)\times_{\Fun(\{1\},\Dcal)}\Ccal $ and $ \Dcal\times \Ccal$ are cartesian if and only if they project to equivalences in $\Ccal$. It follows that $p$ sends $r$-cartesian morphisms to $q$-cartesian morphisms. It remains to show that $p_z$ admits lef adjoint for every $z\in \Dcal.$ Using the equivalence of $\infty$-categories $(\Dcal\times \Ccal)_z\simeq \Ccal$, $p_z$ is the projection 
$$\{z\}\times_{\Fun(\{0\},\Dcal))
}\Fun(\Delta^1,\Dcal)\times_{\Fun(\{1\},\Dcal)}\Ccal\to \Ccal, (z,z\to G(y),y)\mapsto y.$$ By \cite[\href{https://kerodon.net/tag/02J9}{Tag 02J9}]{ker}, it suffices to show that $$(\{z\}\times_{\Fun(\{0\},\Dcal))
}\Fun(\Delta^1,\Dcal)\times_{\Fun(\{1\},\Dcal)}\Ccal)\times_{\Ccal}\times \Ccal_{y/}$$
has an initial object for every $y\in \Ccal.$ By \cite[\href{https://kerodon.net/tag/02VT}{Tag 02VT}]{ker}, there is an equivalence $$\{z\}\times_{\Fun(\{0\},\Dcal))
}\Fun(\Delta^1,\Dcal)\times_{\Fun(\{1\},\Dcal)}\Ccal\simeq \Dcal_{z/}\times_{\Dcal} {\Ccal},$$
so we need to show that $\Ccal_{y/}\times_{\Dcal} \Dcal_{z/}$ has an initial object. Again by \cite[\href{https://kerodon.net/tag/02J9}{Tag 02J9}]{ker}, this is equivalent to showing that for every corepresentable functor $\lambda: \Dcal\to \Scal,$ the composite functor 
$$\Ccal_{y/}\to \Dcal\xrightarrow{\lambda} \Scal$$
is corepresentable. As $\Ccal_{y/}\to \Dcal$ is the composition $\Ccal_{y/}\to \Ccal \xrightarrow{G}\Dcal$ and $G$ admits a left adjoint and so $\lambda \circ G$ is corepresentable, it suffices to show that $\Ccal_{y/}\to \Ccal$ admits a left adjoint. This is true as a left adjoint is given by $y'\mapsto (y\hookrightarrow y\sqcup y').$ Now we know that $p$ admits a relative left adjoint, so the base change of $p$ by $F:\Ecal\to \Dcal$ also admits a left adjoint and we are done. 
\end{proof}
Consequently, we have the following:
\begin{corollary}\label{lm:adjoint_more}
The projection functor $\drpoly_{h,d}(k)\to \Mod_k^\an\times \Alg_k^\an$ admits a left adjoint.
\end{corollary}
We now cite \cite[Corollary 4.7.3.16]{HA} in the following form:
\begin{lemma}\label{lm:adjoint_same}
Let $\Ccal$ be a cocomplete category generated under colimits by $\Ccal^\sfp$, let $\Dcal$ be a cocomplete $\infty$-categories, and let $\Ecal$ be a projectively generated $\infty$-category. Let $G:\Ani(\Ccal)\to \Ecal$ be a functor induced by an inclusion of compact projective generators of $\Ccal$ into the full subcategory of compact projective objects of $\Ecal$. Assume that we have a commutative diagram of $\infty$-categories
$$
\xymatrix{\Ani(\Ccal) \ar[dr]_{G}\ar[r]^-{U} &\Dcal\ar[d]^{G'}\\
&\Ecal
}
$$
satisfying the following conditions:
\begin{enumerate}
    \item $G$ and $G'$ admit left adjoints $F$ and $F',$ respectively.
    \item $G'$ preserves geometric realizations.
    \item $G'$ is conservative.
    \item For each object $x\in S$ a set of compact projective generators of $\Ecal$, the natural morphism $G'F'(x)\to GF(x)$ is an equivalence in $\Ecal.$ 
\end{enumerate}
Then $U$ is an equivalence of $\infty$-categories.
\end{lemma}
\begin{proof}
Let us give some explanation of the assumptions of this lemma. There are $5$ assumptions in \cite[Corollary 4.7.3.16]{HA}:
\begin{enumerate}
    \item This is the same as our assumption $(1)$.
    \item $\Ani(\Ccal)$ admits geometric realizations, so this assumption is always satisfied.
    \item This follows from the assumption that $\Ecal$ is projectively generated and $G'$ preserves geometric realizations.
    \item This assumption is the same as our assumption $(3).$
    \item This is the same as our assumption $(4).$
\end{enumerate}
Moreover, $G$ is fully faithful and hence conservative by the virtue of \cite[Proposition 5.5.8.22]{HTT}.
\end{proof}
\begin{corollary}\label{cor:same}
There exists an equivalence of $\infty$-categories $\drpoly_{h,d}(k)\simeq \Ani(\poly_{h,d}(k)).$
\end{corollary}
\begin{proof}
We have an induced functor $$\Ani(\poly_{h,d}(k))\to \drpoly_{h,d}(k)$$ since each object in a set of compact projective generators $$(K,\T(\Gamma^d(K)))*A,K\hookrightarrow \T(\Gamma^d(K)))*A)$$ for $K\in \poly_{h,d}(k), A\in \Alg_k^\sfp$ is also an object of $\drpoly_{h,d}(k)$ as $\Gamma^d(K)$ continues to be finite projective. We now verify the conditions of Lemma \ref{lm:adjoint_same}:
\begin{enumerate}
    \item This is true by Lemma \ref{lm:adjoint_same} and Lemma \ref{lm:ani_adjoint}.
    \item This is true by Corollary \ref{lm:adjoint_more}.
    \item This is obvious.
    \item This is true since both composite functors map $$(K,A)\mapsto (M,\T(\Gamma^d(M))*A,\Gamma^d(M)\hookrightarrow \T(\Gamma^d(M))*A)\mapsto (M,\T(\Gamma^d(M))*A)$$ for $M\in \Mod_k^\sfp$ and $A\in \Alg_k^\sfp$.\qedhere
\end{enumerate}
\end{proof}
\begin{definition}
We regard the construction that sends each $(B\in \Alg^\an_k)$ to the composite morphism
$$\Lrm\Gamma^d_k(B)\rightarrow \Fr^d(B)\rightarrow B$$
as a functor $\Ffrak_{h,d}: \Alg_k^\an\to \drpoly_{h,d}(k).$ More precisely, it is obtained by animating the ordinary functor $F_{h,d}.$
\end{definition}
\begin{theorem}\label{prop:cliff_main}
The functor $\Ffrak_{h,d}:\Alg_k^\an\to \drpoly_{h,d}(k)$ admits a left adjoint. More explicitly, denote this left adjoint by $\drCl$ and let $\Psi=(K,A,\psi:\Lrm\Gamma^d(K)\to A)\in \drpoly_{h,d}(k)$, then $\drCl(\Psi)$ is given as the pushout of a diagram
$$
\xymatrix{
\Fr(\Lrm\Gamma^d(K)) \ar[r] \ar[d] & A \\
\Fr(K)
}
$$
in the $\infty$-category $\Alg_k^\an.$ Here the morphism $\Fr(\Lrm\Gamma^d(K))\to \Fr(K)$ is induced by the composition $\Lrm\Gamma^d(K)\to \Fr^d(K)\to \Fr(K)$ in $\Mod^\an_k,$ and the morphism $\Fr(\Lrm\Gamma^d(K))\to A$ is induced by $\psi.$
\begin{proof}
By Lemma \ref{lm:ani_adjoint}, Theorem \ref{thm:Cliff}, and Corollary \ref{cor:same}, we obtain the left adjoint $\drCl$ as the animation of functor $\Cl.$ Let $B\in \Alg_k^\an.$ By definition, we have $\Hom_{\Alg_k^\an}(\Cl(\Psi),B)\simeq \Hom_{\Alg_k^\an}(\Psi,\Ffrak_{h,d}(B)).$ By \cite[\href{https://kerodon.net/tag/06BZ}{Tag 06BZ}]{ker}, there is a pullback diagram 
$$
\xymatrix{
\Hom_{\drpoly_{h,d}(k)}(\Psi,\Ffrak_{h,d}(B)) \ar[r] \ar[d] & \Hom_{\Alg_k^\an}(A,B) \ar[d] \\
\Hom_{\Mod^\an_k}(K,B) \ar[r]& \Hom_{\Mod^\an_k}(\Lrm\Gamma^d(K),B)
}
$$
in the $\infty$-category $\Scal.$ As $\Hom_{\Mod^\an_k}(K,B)\simeq \Hom_{\Alg^\an_k}(\Fr(K),B)$ and $\Hom_{\Alg^\an_k}(\Fr(\Lrm\Gamma^d(K)),B),$ this proves the description of $\Cl(\Psi)$ as a pushout in $\Alg^\an_k$.
\end{proof}
\end{theorem}
\begin{proposition}\label{prop:drcClifford}
There is a canonical isomorphism $$\Alg(\pi_0)(\drCl(\Psi))\to \Cl(\pi_0\Psi)$$
for any $\Psi\in \drpoly_{h,d}(k).$
\begin{proof}
As all functors $\Alg(\pi_0), \drCl, \pi_0, \Cl$ preserve colimits, one can assume that $$\Psi=(K,\T(\Gamma^d(K)*A,\Gamma^d(K)\hookrightarrow \T(\Gamma^d(K)*A)$$
for $K\in\Mod_k^\sfp$ and $A\in \Alg_k^\sfp$ so $\pi_0(\Psi)=\Psi$. If $B\in \Alg_k^\sfp$ then $F^{d,\an}(B)\simeq F^d(B)$ so we have 
$$
\begin{matrix*}[l]
\Hom_{\Alg_k}({\Alg(\pi_0)\drCl(\Psi),B})&\simeq& \Hom_{\Alg_k^\an}(\drCl(\Psi),B)\\
&\simeq& \Hom_{\drpoly_{h,d}(k)}(\Psi,F_{h,d}(B))\\
&\simeq&\Hom_{\poly_{h,d}(k)}(\pi_0\Psi,F_{h,d}(B))\\
&\simeq& \Hom_{\poly_{h,d}(k)}(\Cl(\pi_0\Psi),B).
\end{matrix*}
$$
For arbitrary $B\in\Alg_k$ the proof follows by writing $B$ as a colimit of a diagram in $\Alg_k^\sfp$. 
\end{proof}
\end{proposition}
\begin{remark}\label{rmk:higher_homotopy}
\begin{enumerate}
    \item By \cite[Proposition 7.1.4.6]{HA}, we can equip the category of differential graded algebras $\Alg^{\dg}_k$ a model structure where the fibrations are chain maps that are surjective in all degrees and the weak equivalences are quasi-isomorphisms. Some important cofibrations in this model structure are:
    \begin{enumerate}
        \item $k\to k\langle x_n\mid \drm x_n=0\rangle.$ Here $k\langle x\mid \drm x_n=0 \rangle$ denotes the free graded algebra generated by variable $x_n$ in degree $n$ with the differential $\drm$ given by $\drm x_n=0$.
        \item $k\to k\langle x_n,y_{n-1}\mid \drm x_n=y_{n-1}\rangle.$ Here $k\langle x_n,y_{n-1}\mid \drm x_n=y_{n-1}\rangle$ is the free graded algebra generated by $x_n$ in degree $n$ and $y_{n-1}$ in degree $n-1$ with the differential $\drm$ given by $\drm x_{n}=y_{n-1}$.
        \item $k\langle z_{n-1}\mid \drm z_{n-1}=0\rangle\to k\langle x_n,y_{n-1}\mid \drm x_n=y_{n-1}\rangle$ given by $z_{n-1}\mapsto y_{n-1}.$
    \end{enumerate}
    \item In general, derived Clifford algebras are different from ordinary Clifford algebras. We can see such an example in \cite[Example 2.13, Remark 2.14]{Vezz_2016} for derived Clifford algebras of quadratic forms. Based on this, one can give examples for forms of any degree $d\geq 2$. Let $\Psi=0\in \poly_{h,d}(k)$ given by $k\leadsto k, z\mapsto 0$. The derived Clifford algebra of $\Psi$ is given by the pushout diagram
    $$
    \xymatrix{
    k[x] \ar[r] \ar[d] & k \ar[d]\\
    k[x] \ar[r] & \drCl(\Psi)
    }
    $$
in $\Alg^{\an}_k,$ where the upper map is given by $x\mapsto 0$ and the left vertical map is given by $x\mapsto x^d$. One can regard $k[x]$ as $k\langle x=x_0\mid \drm x_0 =0 \rangle.$ Hence $k[x]$ is cofibrant by the above remark and one can compute $\drCl(\Psi)$ by replacing $k[x]\to k$ by its cofibrant replacement and compute the pushout of the new diagram in the ordinary category $\Alg^{\dg}_k.$ Thus we decompose $k[x]\to k$ as follows
$$k[x]\xrightarrow{x\mapsto y} k\langle x_1,y\mid \drm x_1=y\rangle\xrightarrow{x_1,y\mapsto 0} k,$$
in which the first morphism is a cofibration by the above remark and the second morphism is a surjective quasi-isomorphism by direct computations. Then one has the following pushout in $\Alg^{\dg}_k:$
$$
\xymatrix{
k\langle x \rangle \ar[r]^-{x\mapsto y} \ar[d]_{x\mapsto x^d}& k\langle x_1,y\mid \drm x_1 =y \rangle  \ar[d]\\
k\langle x \rangle \ar[r] & \drCl(\Psi).
}
$$
From this description, one sees that \begin{align*}
    \drCl(\Psi)&\simeq k\langle x,x_1,y \mid \drm x_1 = y, \drm x =0, y=x^d  \rangle \\
    &= k\langle x,x_1 \mid \drm x_1 = x^d, \drm x = 0 \rangle.
\end{align*} It has nonzero first homology as $$\drm(x_1x-x x_1)=\drm(x_1)x-x\drm(x_1)=x^{d+1}-x^{d+1}=0,$$ while $x_1 x -x x_1$ is not exact since any element $x^a x_1 x^b x_1 x^c$ in degree $2$ will have differential involving the term $\drm(x_1)=x^d$ which has degree $\geq 2$ in $x$.
\end{enumerate}
\end{remark}
Let $k'$ be a $k$-algebra. As in Section 2, for any algebra $f:k\to k'$ we have a pullback functor $f^*:\poly_{h,d}(k)\to \poly_{h,d}(k')$ and $f_*:\poly_{h,d}(k')\to \poly_{h,d}(k).$ We are going to show there is an adjunction $(f^*\dashv f_*).$ 
\begin{lemma}
There is an adjunction $(f^*\dashv f_*):\Prm_{k'}\rightleftarrows \Prm_{k}.$ Consequently, there exists an adjunction $(\otimes_k k'\dashv (-)_{k}):\poly_{h,d}(k')\to \poly_{h,d}(k).$
\begin{proof}
We already have unit and counit maps $f^*f_* M\to M$ and $N\to f_*f^* N$ for each object $M \in \Prm_{k'}$ and $N\in \Prm_k.$ It remains to show that they are natural. But this is obvious as they are defined by extension of scalars. More precisely, for any $k'$-algebra $A,$ the induced map $$(f^*f_* M)\otimes_{k'} A \to M\otimes_{k'} A$$ is the map $M\otimes_k A\simeq  M\otimes_{k'} (k'\otimes_k A)\to M\otimes_{k'} A$ induced by the ring homomorphism $k'\otimes A\to A, r\otimes a \mapsto ra$; and for any $k$-algebra $A$, the induced map
$$N\otimes_k A \to (f_*f^* N)\otimes_k A$$
is the map $N\otimes_k A \to N\otimes_k (k'\otimes_k A)$ induced by the ring homomorphism $A\to k
\otimes_k A, a\mapsto a\otimes 1.$
\end{proof}
\end{lemma}
It follows that $\otimes_k k'$ preserves colimits and in particular $1$-sifted colimits. This provides a functor $\otimes^{\Lbf}_{k} k'$ $: \drpoly_{h,d}(k)\to \drpoly_{h,d}(k')$ by taking animation of $\otimes_k k'.$ 
\begin{proposition}For any $k$-algebra $k'$, we have $\drCl$$(\Psi\otimes^{\Lbf}_k k')$ $\simeq \drCl$$(\Psi)\otimes^{\Lbf}_k k'.$
\begin{proof}
This is a consequence of \cite[Proposition 5.1.5]{C-S_2019} on animation of composite functors. The only part that we need to prove is $\otimes_R k'$ maps $\poly_{h,d}(k)^{\sfp}$ to $\poly_{h,d}(k')^{\sfp}.$ But this is obvious as $\poly_{h,d}(k)^\sfp$ are retracts of $$\{(M,\T(\Gamma^d(M))*A, \Gamma^d(M)\hookrightarrow T(\Gamma^d(M))*A)\mid M\in \Mod_k^{\sfp}, A\in \Alg_k^{\sfp}\}$$
and $\otimes_R k'$ preserves those objects.
\end{proof}
\end{proposition}
\subsection{Derived Clifford algebras over animated commutative rings}
We now consider the case of any animated commutative ring. Recall from \cite[Section 25.2.1]{spectral} that the $\infty$-category $\scrmod$ is the fiber product $\CAlg^\an\times_{\CAlg(\Sp)}\Mod(\Sp)$ where $\Mod(\Sp)$ is the $\infty$-category of pairs $(R,M)$ where $R$ is an $\Ebb_\infty$-ring and $M$ is an $R$-module. The $\infty$-category $\scrmod^{\mathrm{cn}}$ is defined as the full subcategory of $\scrmod$ spanned by pairs $(R,M)$ where $R\in \CAlg^\an$ and $M$ is a connective $R$-module. If $R$ is an animated commutative ring, we let $\Alg^\an_R$ denote the $\infty$-category $\Alg^\mathrm{cn}_R$ of connective $\Ebb_1$-algebras over $R$ where we regard $R$ as an $\Ebb_2$-ring.

By \cite[Proposition 25.2.1.2]{spectral}, the objects of the full subcategory of $\scrmod^{\mathrm{cn}}$ spanned by pairs $(A,M)$ where $A \simeq \Zbb[x_1,\dots,x_m]$ and $M$ is equivalent to $A^n$ forms compact projective generators of $\scrmod^{\mathrm{cn}}$. Let  $\scrmod^{\cp}$ denote this category. By \cite[Corollary 25.2.1.3]{spectral}, the precomposition with the inclusion $\scrmod^\cp\to \scrmod^\mathrm{cn}$ gives an equivalence
$$\Fun_{\Sigma}(\scrmod^\mathrm{cn},\Ecal)\simeq \Fun(\scrmod^\cp,\Ecal)$$
of $\infty$-categories.
Hence the construction $(A,M)\mapsto \Gamma^d_A(M)$ and $(A,M)\mapsto \T^d_A(M)$ extend to functors $\Lrm\Gamma^d, \Fr^d:\scrmod^\mathrm{cn} \to \scrmod^\mathrm{cn}.$ Moreover, the natural transformation $\Gamma^d_{-}(-) \Rightarrow \T^d_{-}(-)$ induces a natural transformation $\Lrm\Gamma^d_{-}(-)\Rightarrow \Fr^d_{-}(-).$

\begin{definition}
Let $R$ be an animated ring. The $\infty$-category of \textit{derived homogeneous polynomial forms of degree $d$ with noncommutative coefficients}, denoted by $\drpoly_{h,d}(R)$, is the oriented fiber product $\Mod_R^\an \times^{\to}_{\Mod_R^\an} \Alg_R^\an$ where we use functors $\Lrm\Gamma^d_R:\Mod_R^\an\to\Mod_R^\an$ and the forgetful functor $\Alg_R^\an\to \Mod_R^\an$ in the construction.
\end{definition}

\begin{lemma} \label{lm:ani_adj}
Let $\Dcal\times^{\to}_{\Dcal} \Ccal$ be an oriented fiber product with functors $F:\Dcal\to \Dcal$ and $G:\Ccal\to \Dcal$ in its construction. Assume that 
\begin{enumerate}
    \item The $\infty$-categories $\Ccal$ and $\Dcal$ are projectively generated.
    \item Functor $F$ preserves sifted colimits.
    \item $G$ admits a left adjoint and preserves filtered colimits.
\end{enumerate}
Then any functor $H:\Ccal \to \Dcal\times^{\to}_{\Dcal} \Ccal$ satisfying $\ev_0\circ H \simeq G$ and $\ev_1 \circ H \simeq \id_\Ccal$ admits a left adjoint, where $\ev_0: \Dcal\times^{\to}_{\Dcal} \Ccal\to \Dcal$ and $\ev_1: \Dcal\times^{\to}_{\Dcal} \Ccal\to \Ccal$ are projection maps.
\end{lemma}
\begin{proof}
By Lemma \ref{lm:plusdrived_adjunction}, the projection functor $p: \Dcal\times^{\to}_{\Dcal} \Ccal\to \Dcal \times \Ccal$ admits a left adjoint. We now show that $\Dcal\times^{\to}_{\Dcal} \Ccal$ is projectively generated by invoking \cite[Corollary 4.7.3.18]{HA}. First observe that the conditions $(1)$ and $(2)$ imply that $\Dcal\times^{\to}_{\Dcal} \Ccal$ admits filtered colimits and geometric realizations by \cite[\href{https://kerodon.net/tag/06AF}{Tag 06AF}]{ker}. The $\infty$-category $\Dcal\times \Ccal$ is projectively generated by the condition $(1)$. It is also obvious that the projection functor is conservative, and it preserves filtered colimits and geometric realizations by \cite[\href{https://kerodon.net/tag/06AF}{Tag 06AF}]{ker}. 

By Lurie's adjoint functor theorem, it suffices to show that $H$ preserves limits and filtered colimits. Let $q: K^\triangleleft \to \Ccal$ be a limit diagram. Then $\ev_0 \circ H\circ q \simeq G\circ q$ by our assumptions. As $G$ admits a left adjoint, we deduce that $G\circ q$ is a limit diagram. The diagram $\ev_1\circ H \circ q \simeq q$ so it is also a limit diagram. It follows that $H\circ q$ is a limit diagram by \cite[\href{https://kerodon.net/tag/06AF}{Tag 06AF}]{ker}. Similarly, let $q: K^\triangleright \to \Ccal$ be a colimit diagram where $K$ is filtered. Then $\ev_0\circ H \circ q \simeq G\circ q$ so it is a colimit diagram by assumption $(3),$ and $\ev_1\circ G\circ q \simeq q$ is a colimit diagram. Again by \cite[\href{https://kerodon.net/tag/06AF}{Tag 06AF}]{ker}, the diagram $q$ is a colimit diagram and we are done.
\end{proof}

We let $\scralg^\mathrm{cn}$ denote the $\infty$-category of pair $(R,A)$ where $R$ is an animated ring and $A \in \Alg^\an_R$. Let $\scralg^\mathrm{cp}$ denote the full subcategory of $\scralg^\mathrm{cn}$ spanned by pairs $(R,A)$ where $R\simeq \Zbb[x_1,\dots,x_n]$ and $A$ is a finite free $R$-algebra. 
\begin{lemma}\label{lm:scralg}
The objects of $\scralg^\mathrm{cp}$ form compact projective generators for $\scralg^\mathrm{cn}$.
\end{lemma}
\begin{proof}
The proof is similar to that of \cite[Proposition 25.2.1.2]{spectral}. First we show that the objects of  $\scralg^\mathrm{cp}$ are compact projective. As they are finite coproducts of $(\Zbb[x],\Zbb)$ and $(\Zbb,\Zbb\langle x \rangle),$ and compact projective objects are closed under finite coproducts (\cite[Corollary 5.3.4.15, Remark 5.5.8.19]{HTT}), it suffices to show that $(\Zbb[x],\Zbb)$ and $(\Zbb,\Zbb\langle x \rangle)$ are compact projective. This is true as for any $(R,A)\in \scralg^\mathrm{cn},$ 
$$\Hom_{\scralg^\mathrm{cn}}((\Zbb[x],\Zbb),(R,A))\simeq \Hom_{\CAlg^\an}(\Zbb[x],R)\simeq \Omega^\infty(R),$$
where $\Omega^\infty$ sends $R$ to its underlying space in $\Scal$ where we consider $A$ as a connective spectrum. In this connective case, we can also think about $\Omega^\infty$ as the animation of forgetful functors to category $\Set.$ As $\Omega^\infty$ preserves sifted colimits, we deduce that $(\Zbb[x],\Zbb)$ is compact projective. Similarly, 
$$\Hom_{\scralg^\mathrm{cn}}((\Zbb,\Zbb\langle x \rangle),(R,A))\simeq \Hom_{\Alg^\an}(\Zbb\langle x \rangle,A)\simeq \Omega^\infty(A)$$
and hence $(\Zbb,\Zbb\langle x \rangle)$ is compact projective. As the subcategory $\scralg^\cp$ is stable under finite coproducts in $\scralg^\mathrm{cn}$, the inclusion $\scralg^\cp\to \scralg^\mathrm{cn}$ extends to a colimits preserving functor $F:\Pcal_\Sigma(\scralg^\cp)\to \scralg^\mathrm{cn}$ by \cite[Proposition 5.5.8.15]{HTT}, where $\Pcal_\Sigma(\scralg^\cp)$ is the full subcategory of presheaves on $\scralg^\cp$ spanned by those presheaves which preserve finite products. Moreover, the functor $F$ is fully faithful as $\scralg^\cp$ are compact projective. By \cite[Proposition 5.5.8.22]{HTT}, the lemma follows if we can show that $F$ is an equivalence. By Lurie's adjoint functor theorem, the functor $F$ admits a right adjoint $G$. To show that $F$ is an equivalence, it suffices to show that $G$ is conservative. This is true as the functor $$\scralg^\mathrm{cn}\to \Scal\times \Scal, (R,A)\mapsto (\Omega^\infty(R),\Omega^\infty(A))$$
is conservative, and it factors as the composition of its restriction to $\Pcal_\Sigma(\scralg^\cp)$ with $G$. Indeed, 
$$
\begin{matrix*}[l]
(\Omega^\infty(R),\Omega^\infty(A)) & \simeq & (\Hom_{\scralg^\mathrm{cn}}(F(\Zbb,\Zbb\langle x \rangle),(R,A)),\Hom_{\scralg^\mathrm{cn}}(F(\Zbb,\Zbb\langle x \rangle),(R,A))) \\   
&\simeq& (\Hom_{\Pcal_\Sigma(\scralg^\cp)}((\Zbb,\Zbb\langle x \rangle),G(R,A)),\Hom_{\Pcal_\Sigma(\scralg^\cp)}((\Zbb,\Zbb\langle x \rangle),G(R,A)))\\
&\simeq& (\Omega^\infty(G(R,A)),\Omega^\infty(G(R,A))).\hfill\qedhere
\end{matrix*}
$$
\end{proof}

\begin{definition}
We regard the construction that sends each $(B\in \Alg^\an_R)$ to the composite morphism
$$\Lrm\Gamma^d_R(B)\rightarrow \Fr^d_R(B)\rightarrow B$$
as a functor $\Ffrak_{h,d}: \Alg_R^\an\to \drpoly_{h,d}(R).$ More precisely, this is obtained by Lemma \ref{lm:scralg} as the natural transformation $\Gamma^d_{-}(-)\Rightarrow \T^d_{-}(-)\Rightarrow \For$ as functors $\scralg^\cp\to \scrmod^\cp$ induces a natural transformation $\Lrm\Gamma^d_R(-)\Rightarrow \For: \Alg^\an_R\to \Mod^\an_R$ where $\For$ stand for forgetful functors.
\end{definition}

\begin{theorem}
The functor $\frak{F}_{h,d}:\Alg_R^\an\to \drpoly_{h,d}(R)$ admits a left adjoint. More explicitly, denote this left adjoint by $\drCl$ and let $\Psi=(K,A,\psi:\Lrm\Gamma^d(K)\to A)\in \drpoly_{h,d}(R)$. Then $\drCl(\Psi)$ is given as the pushout of a diagram
$$
\xymatrix{
\Fr(\Lrm\Gamma^d(K)) \ar[r] \ar[d] & A \\
\Fr(K)
}
$$
in the $\infty$-category $\Alg_R^\an.$ Here the morphism $\Fr(\Lrm\Gamma^d(K))\to \Fr(K)$ is induced by the composition $\Lrm\Gamma^d(K)\to \Fr^d(K)\to \Fr(K)$ in $\Mod^\an_R,$ and the morphism $\Fr(\Lrm\Gamma^d(K))\to A$ is induced by $\psi.$
\end{theorem}
\begin{proof}
All the assumptions of Lemma \ref{lm:ani_adj} are satisfied, where $F = \Lrm\Gamma^d$ and $G$ is the forgetful functor $\Alg^\an_R\to \Mod_R^\an$. The proof of the description of $\drCl(\Psi)$ as the pushout is similar to that of \ref{prop:cliff_main}.
\end{proof}
\begin{proposition}
For any $R$-algebra $R'$, we have $\drCl(\Psi\otimes^\Lrm_R R')\simeq \drCl(\Psi)\otimes^\Lrm_R R'.$
\begin{proof}
This is a consequence of a more general result that we are going to prove in next section (Proposition \ref{prop:derived_basechange}).
\end{proof}
\end{proposition}
\subsection{Derived Clifford algebras over derived Artin stacks} 
Lurie's straightening functor provides an equivalence
$$\Cart(\CAlg^\an)\xrightarrow{\simeq} \Fun(\CAlg^\an, \Cat_\infty)$$
where $\Cat_\infty$ is the $\infty$-category of $\infty$-categories (\cite[\href{https://kerodon.net/tag/0209}{Tag 0209}]{ker}), and $\Cart(\CAlg^\an)$ is the subcategory of ${\Cat_{\infty}}_{/ \CAlg^\an}$ spanned by cocartesian fibrations over $\CAlg^\an$ and whose morphisms are those functors preserving cocartesian morphisms. Hence by \cite[Proposition 25.2.3.1]{spectral}, the construction $(A,M)\mapsto (A,\Lrm\Gamma(M))$ induces a natural transformation between functors $\Mod^\an_{-}: \CAlg^\an\rightrightarrows \Cat_\infty$. By the universal property of (right) Kan extesions \cite[\href{https://kerodon.net/tag/0309}{Tag 0309}]{ker}, this natural transformation corresponds to a natural transformation between functors $\QCoh(-)^{\mathrm{cn}}: \dst^\alg\rightrightarrows \Cat_\infty$, where $\dst^\alg$ is the $\infty$-category of derived Artin stacks. Roughly speaking, we obtain a collection of functors $$\Lrm\Gamma^d_{\mathfrak{X}}: \QCoh(\Xfrak)^{\mathrm{cn}}\to \QCoh(\Xfrak)^{\mathrm{cn}}$$ 
when derived Artin stacks $\Xfrak$ vary and they are compatible with with pullback functors.
\begin{definition}
Let $\Xfrak$ be a derived Artin stacks. The $\infty$-category of \textit{derived homogeneous polynomial forms of degree $d$ with noncommutative coefficients}, denoted by $\drpoly_{h,d}(\Xfrak)$, is the oriented fiber product $\QCoh(\Xfrak)^{\mathrm{cn}}\times^{\to}_{\QCoh(\Xfrak)^{\mathrm{cn}}} \Alg(\Xfrak)^{\mathrm{cn}}$ where we use functors $\Lrm\Gamma^d_{\mathfrak{X}}: \QCoh(\Xfrak)^{\mathrm{cn}}\to \QCoh(\Xfrak)^{\mathrm{cn}}$ and the forgetful functor $\Alg(\Xfrak)^{\mathrm{cn}}\to \QCoh(\Xfrak)^{\mathrm{cn}}$ in the construction of the oriented fiber product. 
\end{definition}
By \cite[Remark 6.2.2.7]{spectral}, we can identify $\QCoh(\Xfrak)^{\mathrm{cn}}$ with the $\infty$-category $$\Fun^{\ccart}_{\CAlg^\an}(\CAlg(\Xfrak),\scrmod^\mathrm{cn}),$$
the full subcategory of $\Fun_{\CAlg^\an}(\CAlg(\Xfrak),\scrmod^\mathrm{cn})$ spanned by functors $$\CAlg(\Xfrak)\to \scrmod^\mathrm{cn}$$ which maps all morphisms of $\CAlg(\Xfrak)$ to cocartesian morphisms of the cocartesian fibration $\scrmod^\mathrm{cn}\to \CAlg^\an.$ Here $\CAlg(\Xfrak)\to \CAlg^\an$ is the cocartesian fibration corresponding to the functor of points $\CAlg^\an\to \Scal$ of $\Xfrak$. Similarly, we can identify $\Alg(\Xfrak)^\mathrm{cn}$ with $\Fun^{\ccart}_{\CAlg^\an}(\CAlg(\Xfrak),\scralg^\mathrm{cn}).$

As we have an equivalence of $\infty$-categories $$\Fun(\scralg^\mathrm{cn},\Ecal)\to \Fun(\scralg^\cp,\Ecal),$$
the natural transformation given by the construction $(R,A)\mapsto (\Gamma^d(A)\to A^{\otimes d}\to A)$ extends to a natural transformation from $\Lrm\Gamma^d_{-}(-):\scralg^\mathrm{cn}\to \scrmod^\mathrm{cn}$ to the forgetful functor $\scralg^\mathrm{cn}\to \scrmod^\mathrm{cn}.$ This natural transformation corresponds to a functor
$$f_{h,d}:\scralg^\mathrm{cn}\to \scrmod^\mathrm{cn} \times^\to_{\scrmod^\mathrm{cn}} \scralg^\mathrm{cn}$$
where we use $\Lrm\Gamma^d_{-}(-):\scralg^\mathrm{cn}\to \scrmod^\mathrm{cn}$ and the forgetful functor $\scralg^\mathrm{cn}\to \scrmod^\mathrm{cn}$ in the construction of the oriented fiber product.
\begin{definition}
We regard the construction that sends each $(\Bcal\in \Alg(\Xfrak)^\mathrm{cn})$ to the composite morphism
$$\Lrm\Gamma^d_\Xfrak(B)\rightarrow \Fr^d_\Xfrak(B)\rightarrow B$$
as a functor $\Ffrak_{h,d}: \Alg(\Xfrak)^\mathrm{cn}\to \drpoly_{h,d}(\Xfrak).$ More precisely, the functor $f_{h,d}$ gives a commutative square
$$
\xymatrix{
\scralg^\mathrm{cn} \ar[r] \ar[d] & \Fun(\Delta^1,\scrmod^\mathrm{cn}) \ar[d]\\
\scrmod^\mathrm{cn} \times \scralg^{\mathrm{cn}}\ar[r] & \scrmod^\mathrm{cn}\times \scrmod^\mathrm{cn}
}
$$
of simplicial sets, and $\Fun(\CAlg(\Xfrak),-)$ gives the corresponding diagram
$$
\xymatrix{
\Fun^{\ccart}_{\CAlg^\an}(\CAlg(\Xfrak),\scralg^\mathrm{cn}) \ar[d]\ar[r]& \Fun(\Delta^1,\QCoh(\Xfrak)^{\mathrm{cn}})\ar[d] \\
\QCoh(\Xfrak)^\mathrm{cn}\times \Alg(\Xfrak)^\mathrm{cn} \ar[r]& \QCoh(\Xfrak)^\mathrm{cn}\times\QCoh(\Xfrak)^\mathrm{cn}  
}
$$
which induces functor $\Ffrak_{h,d}.$
\end{definition}

\begin{theorem}\label{thm:main_derived}
The functor $\frak{F}_{h,d}:\Alg(\Xfrak)^\mathrm{cn}\to \drpoly_{h,d}(\Xfrak)$ admits a left adjoint. More explicitly, let denote this left adjoint by $\drCl$ and let $\Psi=(\Fcal,\Acal,\psi:\Lrm\Gamma^d(\Fcal)\to \Acal)\in \drpoly_{h,d}(\Xfrak)$. Then $\drCl(\Psi)$ is given as the pushout of a diagram
$$
\xymatrix{
\Fr(\Lrm\Gamma^d(\Fcal)) \ar[r] \ar[d] & \Acal \\
\Fr(\Fcal)
}
$$
in the $\infty$-category $\Alg(\Xfrak)^{\mathrm{cn}}.$ Here the morphism $\Fr(\Lrm\Gamma^d(\Fcal))\to \Fr(\Fcal)$ is induced by the composition $\Lrm\Gamma^d(\Fcal)\to \Fr^d(\Fcal)\to \Fr(\Fcal)$ in $\Mod^\an_R,$ and the morphism $\Fr(\Lrm\Gamma^d(\Fcal))\to \Acal$ is induced by $\psi.$
\end{theorem}
\begin{proof}
By \cite[\href{https://kerodon.net/tag/02J9}{Tag 02J9}]{ker}, it suffices to show that for any $\Psi\in \drpoly_{h,d}(\Xfrak),$ the composite functor $\Hom_{\drpoly_{h,d}(\Xfrak)}(\Psi,\Ffrak_{h,d}(-))$ is corepresentable. By \cite[\href{https://kerodon.net/tag/06BZ}{Tag 06BZ}]{ker}, there is a pullback diagram
$$
\xymatrix{
\Hom_{\drpoly_{h,d}(\Xfrak)}(\Psi,\Ffrak_{h,d}(-)) \ar[r] \ar[d] & \Hom_{\Alg(\Xfrak)^\mathrm{cn}}(\Acal,\ev_1(\Ffrak_{h,d}(-))) \ar[d] \\
\Hom_{\QCoh(\Xfrak)^\mathrm{cn}}(\Fcal,\ev_0(\Ffrak_{h,d}(-))) \ar[r]& \Hom_{\QCoh(\Xfrak)^{\mathrm{cn}}}(\Lrm\Gamma^d(\Fcal),\ev_1(\Ffrak_{h,d}(-)))
}
$$
in the $\infty$-category $\Fun(\Alg(\Xfrak)^\mathrm{cn},\Scal).$ As $\ev_{i}(\Ffrak_{h,d}(\Bcal))=\Bcal$, we obtain the pullback diagram
$$
\xymatrix{
\Hom_{\drpoly_{h,d}(\Xfrak)}(\Psi,\Ffrak_{h,d}(-)) \ar[r] \ar[d] & \Hom_{\Alg(\Xfrak)^\mathrm{cn}}(\Acal,-) \ar[d] \\
\Hom_{\QCoh(\Xfrak)^\mathrm{cn}}(\Fcal,-) \ar[r]& \Hom_{\QCoh(\Xfrak)^{\mathrm{cn}}}(\Lrm\Gamma^d(\Fcal),-).
}
$$
As the construction $\Bcal\mapsto (\Lrm\Gamma^d(\Fcal)\to \Acal\to \Bcal)$ provides the morphism $$\Hom_{\Alg(\Xfrak)^\mathrm{cn}}(\Acal,-) \to \Hom_{\QCoh(\Xfrak)^\mathrm{cn}}(\Fcal,-),$$ it is equivalence to the morphism 
$$\Hom_{\Alg(\Xfrak)^\mathrm{cn}}(\Acal,-) \to \Hom_{\Alg(\Xfrak)^\mathrm{cn}}(\Fr(\Fcal),-)$$
which is induced by the composition with the morphism $\Fr(\Fcal)\to \Acal$ in $\Alg(\Xfrak)^\mathrm{cn}.$ On the other hand, the morphism 
$$\Hom_{\QCoh(\Xfrak)^\mathrm{cn}}(\Fcal,-) \to \Hom_{\QCoh(\Xfrak)^{\mathrm{cn}}}(\Lrm\Gamma^d(\Fcal),-)$$ 
is equivalent to the morphism 
$$\Hom_{\Alg(\Xfrak)^\mathrm{cn}}(\Fr(\Fcal),-) \to \Hom_{\Alg(\Xfrak)^{\mathrm{cn}}}(\Fr(\Lrm\Gamma^d(\Fcal)),-)$$
which is induced by the composition with the morphism $\Fr(\Lrm\Gamma^d(\Fcal))\to \Fr(\Fcal)$ given in our statement. Therefore we obtain the pullback diagram
$$
\xymatrix{
\Hom_{\drpoly_{h,d}(\Xfrak)}(\Psi,\Ffrak_{h,d}(-)) \ar[r] \ar[d] & \Hom_{\Alg(\Xfrak)^\mathrm{cn}}(\Acal,-) \ar[d] \\
\Hom_{\Alg(\Xfrak)^\mathrm{cn}}(\Fr(\Fcal),-) \ar[r]& \Hom_{\Alg(\Xfrak)^{\mathrm{cn}}}(\Fr(\Lrm\Gamma^d(\Fcal)),-).
}
$$
Thus, the functor $\Hom_{\drpoly_{h,d}(\Xfrak)}(\Psi,\Ffrak_{h,d}(-))$ is corepresented by the pushout of the corresponding diagram 
$$
\xymatrix{
\Fr(\Lrm\Gamma^d(\Fcal)) \ar[r] \ar[d] & \Acal \\
\Fr(\Fcal).
}
$$
\end{proof}
\begin{proposition}\label{prop:derived_basechange}
For any morphism $f:\Xfrak\to \Yfrak$ of derived Artin stacks, we have an isomorphism $f^*\drCl(\Psi)\simeq \drCl(f^*\Psi)$ for $\Psi\in \drpoly_{h,d}(\Yfrak).$
\begin{proof}
Pick any $\Psi=(\Fcal,\Acal,\psi:\Lrm\Gamma^d(\Fcal)\to \Acal)\in \drpoly_{h,d}(\Yfrak).$ By Theorem \ref{thm:main_derived}, we have the pushout diagram
$$
\xymatrix{
\Fr(\Lrm\Gamma^d(\Fcal)) \ar[r] \ar[d] & \Acal\ar[d] \\
\Fr(\Fcal) \ar[r]& \drCl(\Psi).
}
$$
Since $f^*$ commute with colimits, we obtain the pushout diagram
$$
\xymatrix{
\Fr(\Lrm\Gamma^d(f^*\Fcal)) \ar[r] \ar[d] & f^*\Acal\ar[d] \\
\Fr(f^*\Fcal) \ar[r]& f^*\drCl(\Psi)
}
$$
where $\Fr(\Lrm\Gamma^d(f^*\Fcal))\to f^*\Acal$ is induced by $f^*\psi$. Thus this pushout diagram computes $\Cl(f^*\Psi)$ and we are done.
\end{proof}
\end{proposition}
\bibliographystyle{style}

\end{document}